\documentclass[12pt]{amsart}
\input{epsf}
\usepackage{amssymb,latexsym,cite}
\topskip  \numberwithin{equation}{section}
\newtheorem{theorem}{Theorem}[section]
\newtheorem{proposition}[theorem]{Proposition}

\newtheorem{remark}[theorem]{Remark}

\begin{document}

\pagenumbering{arabic}
\pagestyle{headings}

\newcommand{\DPB}[4]{P\beta_{#1}^{(#2)}(#3,#4)}

\title[Stirling numbers associated with sequences of polynomials]{Stirling numbers associated with sequences of polynomials}
\author{Dae San Kim}
\address{Department of Mathematics, Sogang University, Seoul 121-742, Republic of Korea}
\email{dskim@sogang.ac.kr}

\author{Taekyun Kim}
\address{Department of Mathematics, Kwangwoon University, Seoul 139-701, Republic of Korea}
\email{tkkim@kw.ac.kr}


\begin{abstract}
Let $\bold{P}=\left\{p_{n}(x)\right\}_{n=0}^{\infty}$ be a sequence of polynomials such that deg\,$p_{n}(x)=n,\, p_{0}(x)=1$. The aim of this paper is to study the Stirling numbers of the second kind associated with $\bold{P}$ and of the first kind associated with $\bold{P}$, in a unified and systematic way with the help of umbral calculus technique. This generalizes enormously the `classical' Stirling numbers of both kinds, which correspond to the sequence $\left\{x^{n}\right\}_{n=0}^{\infty}$. Our results are illustrated with many examples which give rise to interesting inverse relations in each case.
\end{abstract}
\maketitle

\markboth{\centerline{\scriptsize Stirling numbers associated with sequences of polynomials}}{\centerline{\scriptsize Dae San Kim, Taekyun Kim}}


\section{Introduction}

In his famous paper [6] of 1749,  Euler discovered a method of computing values of the Riemann zeta function $\zeta(s)$ at negative integers, although he never found a convergent analytic expression for $\zeta(s)$ at negative integers.
In fact, he introduced the Euler zeta function $\eta(s)$ given by
\begin{align*}
\eta(s)=1-\frac{1}{2^{s}}+\frac{1}{3^{s}}-\frac{1}{4^{s}}+\cdots.   
\end{align*}
This series converges for $\mathrm{Re}(s) > 0$. The function $\eta(s)$ can be analytically continued to the whole complex plane. It is related to the $\zeta(s)$ by $\eta(s)=(1-2^{1-s})\zeta(s)$, for $\mathrm{Re}(s) >1$. For the purpose of computing $\eta(-n)$, for $n=0,1,2,\dots$, he was led to introduce the Eulerian polynomials $A_{n}(t)$ by
\begin{align*}
\frac{A_{n}(x)}{(1-x)^{n+1}}=\sum_{j=0}^{\infty}x^{j}(j+1)^{n},
\end{align*}
the coefficients of which are the Eulerian numbers.
Eulerian numbers and polynomials have long been studied due to their number-theoretic and combinatorial importance, since Euler's introduction. There are many ways of defining the Eulerian polynomials and numbers,  including the ones given in (a)-(h) and (k) in Proposition 3.2. Let $[n]=\left\{1,2,\dots,n\right\}$. A combinatorial interpretation of the Eulerian polynomial $A_{n}(t)$ is as follows. For each permutation $\sigma$ in the symmetric group $S_{n}$, the descent and the excedance of $\sigma$ are respectively given by
\begin{align*}
&D(\sigma)= \left\{i \in [n-1] : \sigma(i) > \sigma(i + 1)\right\}, \\
&E(\sigma)= \left\{i \in [n-1] : \sigma(i) > i\right\}.
\end{align*}
Then $A_{n}(t)$ is given by 
\begin{align*}
A_{n}(t)=\sum_{\sigma \in S_{n}} t^{d(\sigma)}=\sum_{\sigma \in S_{n}} t^{e(\sigma)},
\end{align*}
where $d(\sigma)=|D(\sigma)|, e(\sigma)=|E(\sigma)|$ (see [2,7,8,21,25]).\par
The Stirling number of the second $S_{2}(n,k)$ is the number of ways to partition a set of $n$ objects into $k$ nonempty subsets. Eulerian numbers and the Stirling numbers of the second kind are closely related by the relation due to Frobenius (see Proposition 3.2, (k), (l)).

Let $\bold{P}=\left\{p_{n}(x)\right\}_{n=0}^{\infty}$ be a sequence of polynomials such that deg\,$p_{n}(x)=n, p_{0}(x)=1$.
In [19], Koutras generalized the classical Eulerian numbers and polynomials to $p_{n}$-associated Eulerian numbers and polynomials with motivation of providing a unified approach to the study of Eulerian-related numbers and with combinatorial, probabilistic and statistical applications in mind. In addition, he mentioned in the same paper that a lot of special numbers, like the Stirling numbers of the second, the Lah numbers and Gould-Hopper numbers, appear in the expansions of some polynomials in terms of the falling factorials. This observation is the motivation and impetus for the present research on the Stirling numbers of both kinds associated with sequences of polynomials (see [19,20]). \\
The aim of this paper is to study the Stirling numbers of the second $S_{2}(n,k;\bold{P})$ associated with any sequence of polynomials $\bold{P}$, which are defined as the coefficients in the expansion of $p_{n}(x)$ in terms of the falling factorials $(x)_{k}$. This is done with the help of umbral calculus technique. We illustrate our results with twenty examples. Moreover, we also introduce the Stirling numbers of the first kind $S_{1}(n,k;\bold{P})$ associated with any sequence of polynomials $\bold{P}$, which are defined as the coefficients in the expansion of $(x)_{n}$ in terms of the polynomials $p_{k}(x)$.  Again, we study them by means of umbral calculus and illustrate our results with the same twenty examples for the Striling numbers of the second kind.
We note here that the inverse relations for the ordinary Stirling numbers also hold for the Stirling numbers associated with sequences of polynomials $\bold{P}$. This observation applied to each individual case yields interesting orthogonality and inverse relations for the twenty examples, which are explicitly stated. For example, in the case of the sequence  $\bold{P}=\left\{B_{n}(x)\right\}$ of Bernoulli polynomials, we obtain:
\begin{align*}
&\sum_{k=l}^{n}\sum_{m=k}^{n}\sum_{j=l}^{k}\frac{1}{k!}\binom{k}{j}\frac{(m+1)_{k}}{m+1}S_{1}(n,m)S_{2}(j,l)B_{k-j}=\delta_{n,l},\\
&\sum_{k=l}^{n}\sum_{j=k}^{n}\sum_{m=l}^{k}\frac{1}{l!}\binom{n}{j}\frac{(m+1)_{l}}{m+1}S_{2}(j,k)S_{1}(k,m)B_{n-j}=\delta_{n,l},\\
&a_{n}=\sum_{k=0}^{n}\sum_{l=k}^{n}S_{2}(l,k)\binom{n}{l}B_{n-l}c_{k} \Longleftrightarrow c_{n}=\sum_{k=0}^{n}\frac{1}{k!}\sum_{l=k}^{n}\frac{S_{1}(n,l)}{l+1}(l+1)_{k}a_{k},\\
&a_{n}=\sum_{k=n}^{m}\sum_{l=n}^{k}S_{2}(l,n)\binom{k}{l}B_{k-l}c_{k} \Longleftrightarrow c_{n}=\sum_{k=n}^{m}\frac{1}{n!}\sum_{l=n}^{k}\frac{S_{1}(k,l)}{l+1}(l+1)_{n}a_{k},
\end{align*}
where $S_{1}(n,k), S_{2}(n,k), B_{n}, (x)_{n}$ are the Stirling numbers of the first kind, the Stirling numbers of the second kind, Bernoulli numbers and the falling factorials, respectively.
As another example, in the case of the sequence $\bold{P}=\left\{x^{[n]}\right\}$ of central factorials we have:
\begin{align*}
&\sum_{k=l}^{n}\sum_{m=k}^{n}\sum_{j=l}^{k}S_{1}(n,m)T_{2}(m,k)T_{1}(k,j)S_{2}(j,l)=\delta_{n,l}, \\
&\sum_{k=l}^{n}\sum_{j=k}^{n}\sum_{m=l}^{k}T_{1}(n,j)S_{2}(j,k)S_{1}(k,m)T_{2}(m,l)=\delta_{n,l},\\
&a_{n}=\sum_{k=0}^{n}\sum_{l=k}^{n}T_{1}(n,l)S_{2}(l,k)c_{k} \Longleftrightarrow c_{n}=\sum_{k=0}^{n}\sum_{l=k}^{n}S_{1}(n,l)T_{2}(l,k)a_{k},\\
&a_{n}=\sum_{k=n}^{m}\sum_{l=n}^{k}T_{1}(k,l)S_{2}(l,n)c_{k} \Longleftrightarrow c_{n}=\sum_{k=n}^{m}\sum_{l=n}^{k}S_{1}(k,l)T_{2}(l,n)a_{k},
\end{align*}
where $T_{1}(n,k)$ and $T_{2}(n,k)$ are the central factorial numbers of the first kind and the central factorial numbers of the second kind.
In this investigation of Stirling numbers, we are led to discover natural definitions for logarithms and exponentials associated to any delta series which become even clearer when we study the central factorial numbers of both kinds associated with $\bold{P}$ in our forthcoming paper. In addition, we also go over the Eulerian numbers and polynomials associated with $\bold{P}$, of which some are not considered by Koutras [19]. 
The novelty of this paper is that it is the first paper which studies the Stirling numbers of both kinds associated with any sequence of polynomials in a unified and systematic way with the help of umbral calculus. \\
The outline of this paper is as follows. In Section 2, we briefly go over umbral calculus. In Section 3, we introduce the Stirling numbers of the second kind associated with sequences of polynomials. We illustrate our results with twenty examples in Section 4. In Section 5, we introduce the Stirling numbers of the first kind associated with sequences of polynomials. We illustrate our results in Section 6 with the same examples from Section 4 and thereby get interesting inverse relations for each one. Finally, we conclude our paper in Section 7.

\section{Review of umbral calculus}
\vspace{0.5cm}
Here we will briefly go over very basic facts about umbral calculus. For more details on this, we recommend the reader to refer to [1,5,23,24]. We remark that recently umbral calculus has been extended to the case of degenerate umbral calculus in order to treat degenerate special polynomials and numbers, which involves degenerate exponentials (see [12]).
Let $\mathbb{C}$ be the field of complex numbers. Then $\mathcal{F}$ denotes the algebra of formal power series in $t$ over $\mathbb{C}$, given by
\begin{displaymath}
 \mathcal{F}=\bigg\{f(t)=\sum_{k=0}^{\infty}a_{k}\frac{t^{k}}{k!}~\bigg|~a_{k}\in\mathbb{C}\bigg\},
\end{displaymath}
and $\mathbb{P}=\mathbb{C}[x]$ indicates the algebra of polynomials in $x$ with coefficients in $\mathbb{C}$. \par
Let $\mathbb{P}^{*}$ be the vector space of all linear functionals on $\mathbb{P}$. If $\langle L|p(x)\rangle$ denotes the action of the linear functional $L$ on the polynomial $p(x)$, then the vector space operations on $\mathbb{P}^{*}$ are defined by
\begin{displaymath}
\langle L+M|p(x)\rangle=\langle L|p(x)\rangle+\langle M|p(x)\rangle,\quad\langle cL|p(x)\rangle=c\langle L|p(x)\rangle,
\end{displaymath}
where $c$ is a complex number. \par
For $f(t)\in\mathcal{F}$ with $\displaystyle f(t)=\sum_{k=0}^{\infty}a_{k}\frac{t^{k}}{k!}\displaystyle$, we define the linear functional on $\mathbb{P}$ by
\begin{equation}\label{1B}
\langle f(t)|x^{k}\rangle=a_{k}. 
\end{equation}
From \eqref{1B}, we note that
\begin{equation*}
 \langle t^{k}|x^{n}\rangle=n!\delta_{n,k},\quad(n,k\ge 0), 
\end{equation*}
where $\delta_{n,k}$ is the Kronecker's symbol. \par
Some remarkable linear functionals are as follows:
\begin{align}
&\langle e^{yt}|p(x) \rangle=p(y), \nonumber \\
&\langle e^{yt}-1|p(x) \rangle=p(y)-p(0), \label{2B} \\
& \bigg\langle \frac{e^{yt}-1}{t}\bigg |p(x) \bigg\rangle = \int_{0}^{y}p(u) du.\nonumber
\end{align}
Let
\begin{equation}\label{3B}
 f_{L}(t)=\sum_{k=0}^{\infty}\langle L|x^{k}\rangle\frac{t^{k}}{k!}.
\end{equation}
Then, by \eqref{1B} and \eqref{3B}, we get
\begin{displaymath}
 \langle f_{L}(t)|x^{n}\rangle=\langle L|x^{n}\rangle.
\end{displaymath}
That is, $f_{L}(t)=L$, as linear functionals on $\mathbb{P}$. In fact, the map $L\longmapsto f_{L}(t)$ is a vector space isomorphism from $\mathbb{P}^{*}$ onto $\mathcal{F}$.\par  Henceforth, $\mathcal{F}$ denotes both the algebra of formal power series  in $t$ and the vector space of all linear functionals on $\mathbb{P}$. $\mathcal{F}$ is called the umbral algebra and the umbral calculus is the study of umbral algebra. 
For each nonnegative integer $k$, the differential operator $t^k$ on $\mathbb{P}$ is defined by
\begin{equation}\label{4B}
t^{k}x^n=\left\{\begin{array}{cc}
(n)_{k}x^{n-k}, & \textrm{if $k\le n$,}\\
0, & \textrm{if $k>n$.}
\end{array}\right. 
\end{equation}
Extending \eqref{4B} linearly, any power series
\begin{displaymath}
 f(t)=\sum_{k=0}^{\infty}\frac{a_{k}}{k!}t^{k}\in\mathcal{F}
\end{displaymath}
gives the differential operator on $\mathbb{P}$ defined by
\begin{equation}\label{5B}
 f(t)x^n=\sum_{k=0}^{n}\binom{n}{k}a_{k}x^{n-k},\quad(n\ge 0). 
\end{equation}
It should be observed that, for any formal power series $f(t)$ and any polynomial $p(x)$, we have
\begin{equation}\label{6B}
\langle f(t) | p(x) \rangle =\langle 1 | f(t)p(x) \rangle =f(t)p(x)|_{x=0}.
\end{equation}
Here we note that an element $f(t)$ of $\mathcal{F}$ is a formal power series, a linear functional and a differential  operator. Some notable differential operators are as follows: 
\begin{align}
&e^{yt}p(x)=p(x+y), \nonumber\\
&(e^{yt}-1)p(x)=p(x+y)-p(x), \label{7B}\\
&\frac{e^{yt}-1}{t}p(x)=\int_{x}^{x+y}p(u) du.\nonumber
\end{align}

The order $o(f(t))$ of the power series $f(t)(\ne 0)$ is the smallest integer for which $a_{k}$ does not vanish. If $o(f(t))=0$, then $f(t)$ is called an invertible series. If $o(f(t))=1$, then $f(t)$ is called a delta series. \par
For $f(t),g(t)\in\mathcal{F}$ with $o(f(t))=1$ and $o(g(t))=0$, there exists a unique sequence $s_{n}(x)$ (deg\,$s_{n}(x)=n$) of polynomials such that
\begin{equation} \label{8B}
\big\langle g(t)f(t)^{k}|s_{n}(x)\big\rangle=n!\delta_{n,k},\quad(n,k\ge 0).
\end{equation}
The sequence $s_{n}(x)$ is said to be the Sheffer sequence for $(g(t),f(t))$, which is denoted by $s_{n}(x)\sim (g(t),f(t))$. We observe from \eqref{8B} that 
\begin{equation}\label{9B}
s_{n}(x)=\frac{1}{g(t)}q_{n}(x),
\end{equation}
where $q_{n}(x)=g(t)s_{n}(x) \sim (1,f(t))$.\par
In particular, if $s_{n}(x) \sim (g(t),t)$, then $q_{n}(x)=x^n$, and hence 
\begin{equation}\label{10B}
s_{n}(x)=\frac{1}{g(t)}x^n.
\end{equation}

It is well known that $s_{n}(x)\sim (g(t),f(t))$ if and only if
\begin{equation}\label{11B}
\frac{1}{g\big(\overline{f}(t)\big)}e^{x\overline{f}(t)}=\sum_{k=0}^{\infty}\frac{s_{k}(x)}{k!}t^{k}, 
\end{equation}
for all $x\in\mathbb{C}$, where $\overline{f}(t)$ is the compositional inverse of $f(t)$ such that $\overline{f}(f(t))=f(\overline{f}(t))=t$. \par

The following equations \eqref{12B}, \eqref{13B}, and \eqref{14B} are equivalent to the fact that  $s_{n}\left(x\right)$ is Sheffer for $\left(g\left(t\right),f\left(t\right)\right)$, for some invertible $g(t)$: 
\begin{align}
f\left(t\right)s_{n}\left(x\right)&=ns_{n-1}\left(x\right),\quad\left(n\ge0\right),\label{12B}\\
s_{n}\left(x+y\right)&=\sum_{j=0}^{n}\binom{n}{j}s_{j}\left(x\right)q_{n-j}\left(y\right),\label{13B}
\end{align}
with $q_{n}\left(x\right)=g\left(t\right)s_{n}\left(x\right),$
\begin{equation}\label{14B}
s_{n}\left(x\right)=\sum_{j=0}^{n}\frac{1}{j!}\big\langle{g\left(\overline{f}\left(t\right)\right)^{-1}\overline{f}\left(t\right)^{j}}\big |{x^{n}\big\rangle}x^{j}.
\end{equation}
If $s_{n}(x)\sim(g(t),f(t))$, then the following recurrence relation holds:
\begin{equation}\label{15B}
s_{n+1}(x)=\bigg(x-\frac{g'(t)}{g(t)}\bigg)\frac{1}{f'(t)}s_{n}(x).
\end{equation}
Let $s_n(x) \sim (g(t),f(t))$. Then, for any polynomial $p(x)$, we have the polynomial exapnsion given by
\begin{align}
p(x)=\sum_{k \ge 0}\frac{1}{k!}\big\langle{g(t)f(t)^{k}|p(x)\big\rangle}s_{k}(x).\label{16B}
\end{align}
For any formal power series $h(t)$, any polynomial $p(x)$ and $ a \in \mathbb{C}$, we have
\begin{align}
\langle{h(at) | p(x) \rangle}=\langle{h(t) | p(ax) \rangle}.\label{17B}
\end{align}

\section{Stirling numbers of the second kind associated with sequences of polynomials}
We recall that the Stirling numbers of the second kind $S_2(n,k)$ are given by 
\begin{align}
x^{n}=\sum_{k=0}^{n}S_2(n,k)(x)_k, \label{1C}
\end{align}
where the falling factorials $(x)_{n}$ are defined by
\begin{align}
(x)_{n}=x(x-1)\cdots(x-n+1), \,\,(n \ge 1),\quad(x)_{0}=1.\label{2C}
\end{align}
More generally, for any real number $\lambda$, the generalized falling factorials are defined by
\begin{align}
(x)_{n,\lambda}=x(x-\lambda)\cdots(x-(n-1)\lambda), \,\,(n \ge 1),\quad(x)_{0,\lambda}=1.\label{3C}
\end{align}
Then, as a degenerate version of the Stirling numbers of the second kind, the degenerate Stirling numbers of the second kind $S_{2,\lambda}(n,k)$ are given by
\begin{align}
(x)_{n,\lambda}=\sum_{k=0}^{n}S_{2,\lambda}(n,k)(x)_{k}.\label{4C}
\end{align}
Let $\bold{P}=\left\{p_{n}(x)\right\}_{n=0}^{\infty}$ be a sequence of polynomials such that deg\,$p_{n}(x)=n, p_{0}(x)=1$.
In \eqref{1C} and \eqref{4C}, we observe that $S_{2}(n,k)$ and $S_{2,\lambda}(n,k)$ arise as the coefficients respectively when we expand $x^{n}$ and $(x)_{n,\lambda}$ in terms of $(x)_{k}$. In addition, the (unsigned) Lah number $L(n,k)$ and Gould-Hopper number $G(n,k;r,s)$ also appear as the coefficients in the expansion of the rising factorial $\langle{x\rangle}_{n}$ and $(rx+s)_{n}$ in terms of $(x)_{k}$, where $r \ne 0$. In view of this observation, it seems natural to define the {\it{Stirling numbers of the second associated with $\bold{P}=\left\{p_{n}(x)\right\}_{n=0}^{\infty}$}} as the coefficients when we expand $p_{n}(x)$ in terms of $(x)_{k}$:
\begin{align}
p_{n}(x)=\sum_{k=0}^{n}S_2(n,k;\bold{P})(x)_{k}.\label{5C}
\end{align}

\begin{theorem}
Let $\bold{P}=\left\{p_{n}(x)\right\}_{n=0}^{\infty}$ be a sequence of polynomials such that deg\,$p_{n}(x)=n, p_{0}(x)=1$, with $p_{n}(x)=\sum_{l=0}^{n}p_{n,l}x^{l}$. We let
\begin{align*}
p_{n}(x)=\sum_{k=0}^{n}S_2(n,k;\bold{P})(x)_{k}.
\end{align*}
(a)\,\, Then the Stirling numbers of the second kind associated with $\bold{P}$ are given by
\begin{align*}
S_2(n,k;\bold{P})=\frac{1}{k!}\langle{(e^{t}-1)^{k}|p_{n}(x)\rangle}.
\end{align*}
More explicitly, it is given by
\begin{align*}
S_2(n,k;\bold{P})=\sum_{l=k}^{n}S_{2}(l,k)p_{n,l}.
\end{align*}
(b)\,\, $p_{n,l}=\sum_{k=l}^{n}S_{1}(k,l)S_{2}(n,k;\bold{P})$. \\
(c)\,\,Let $p_{n}(x)$ be Sheffer for the pair $(g(t), f(t))$. Then the generating function of $S_2(n,k;\bold{P})$ is given by
\begin{align*}
\sum_{n=k}^{\infty}S_2(n,k;\bold{P})\frac{t^{n}}{n!}=\frac{1}{g(\bar{f}(t))}\frac{1}{k!}(e^{\bar{f}(t)}-1)^{k}=\frac{1}{g(\bar{f}(t))}\sum_{n=k}^{\infty}S_{2}(n,k)\frac{\bar{f}(t)^{n}}{n!}.
\end{align*}
(d),\, \textrm{Let $\bar{\bold{P}}=\left\{\bar{p}_{n}(x)\right\}$, where $\bar{p}_{n}(x)$ is the sequence of polynomials defined by} $\bar{p}_{0}(x)=1,\,\,\bar{p}_{n}(x)=xp_{n-1}(x),\,\,(n \ge 1)$. Then we have
\begin{align*}
S_{2}(n+1,k;\bar{\bold{P}})=S_{2}(n,k-1;\bold{P})+kS_{2}(n,k;\bold{P}),\,\,(0 \le k \le n+1).
\end{align*}
\begin{proof}
(a)\,\, By noting $(x)_{n} \sim (1, e^{t}-1)$, this follows from \eqref{16B}. \\
From \eqref{16B} and \eqref{5C}, we have
\begin{align*}
S_{2}(n,k;\bold{P})&=\frac{1}{k!}\big\langle{(e^{t}-1)^{k}|p_{n}(x)\big\rangle}\\
&=\sum_{l=k}^{n}S_{2}(l,k)\frac{1}{l!}\big(\frac{d}{dx}\big)^{l}p_{n}(x)|_{x=0}\\
&=\sum_{l=k}^{n}S_{2}(l,k)p_{n,l}.
\end{align*}
(b) The identity is immediate from the following observation:
\begin{align*}
\sum_{l=0}^{n}p_{n,l}x^{l}&=\sum_{k=0}^{n}S_{2}(n,k;\bold{P})(x)_{k}\\
&=\sum_{k=0}^{n}S_{2}(n,k;\bold{P})\sum_{l=0}^{k}S_{1}(k,l)x^{l}\\
&=\sum_{l=0}^{n}\sum_{k=l}^{n}S_{1}(k,l)S_{2}(n,k;\bold{P})x^{l}.
\end{align*}
Note here that the explicit expression of $S_{2}(n,k;\bold{P})$ also follows from this by inversion.
(c)\,\, From \eqref{11B} and \eqref{5C}, we have
\begin{align*}
\sum_{k=0}^{\infty}\sum_{n=k}^{\infty}S_2(n,k;\bold{P})\frac{t^{n}}{n!}(u)_{k}&=\sum_{n=0}^{\infty}\sum_{k=0}^{n}S_{2}(n,k;\bold{P})(u)_{k}\frac{t^{n}}{n!}\\
&=\sum_{n=0}^{\infty}p_{n}(u)\frac{t^{n}}{n!}=\frac{1}{g(\bar{f}(t))}\big(1+(e^{\bar{f}(t)}-1)\big)^{u}\\
&=\sum_{k=0}^{\infty}\frac{1}{g(\bar{f}(t))}\frac{1}{k!}(e^{\bar{f}(t)}-1)^{k}(u)_{k}.
\end{align*}
(d)\,\, This is immediate from the following observation:
\begin{align*}
\sum_{k=0}^{n+1}S_{2}(n+1,k;\bar{\bold{P}})(x)_{k}&=xp_{n}(x)=
\sum_{k=0}^{n}(x-k+k)S_{2}(n,k;\bold{P})(x)_{k}\\
&=\sum_{k=0}^{n}S_{2}(n,k;\bold{P})(x)_{k+1}+\sum_{k=0}^{n}kS_{2}(n,k;\bold{P})(x)_{k}\\
&=\sum_{k=0}^{n+1}S_{2}(n,k-1;\bold{P})(x)_{k}+\sum_{k=0}^{n+1}kS_{2}(n,k;\bold{P})(x)_{k}.
\end{align*} 
\end{proof}
\end{theorem}

Here we recall some properties of the `classical' Eulerian numbers in Proposition 3.2,  which correspond to the sequence $\left\{x^{n}\right\}_{n=0}^{\infty}$, and their generalization to any sequence of polynomials by Koutras (see [19]). The Stirling numbers of the second kind associated with $\bold{P}=\left\{p_{n}(x)\right\}$ are related in a simple way to the Eulerian numbers associated with $\bold{P}=\left\{p_{n}(x)\right\}$ (see Theorem 3.3, (h)). This generalizes the relation between Stirling numbers of the second kind and Eulerian numbers which is due to Frobenius (see Proposition 3.2, (k)).\\
The classical Eulerian polynomial $A_{n}(x)=\sum_{k=0}^{n-1}A_{n,k}x^{k}$ is monic of degree $n-1$, for $n \ge 1$, with $A_{0}(x)=1$. They are recursively defined by the relation in (a) below. As a matter of convenience, we write $A_{n}(x)=\sum_{k=0}^{n}A_{n,k}x^{k},\,\,(n \ge 0)$, with the understanding that $A_{n,n}=0$, for $n \ge 1$, and $A_{0,0}=1$. For more details on classical Eulerian numbers and polynomials, we let the reader refer to [7], while for their generalizations to any sequence of polynomials we recommend the reader to refer to the paper [19]. Here the reader ought to note that the definition of Eulerian polynomials associated with sequence of polynomials in [19] are different from ours. Indeed, when $\bold{P}=\left\{x^{n}\right\}$, namely in the case of classical Eulerian polynomials, his definition of the $n$th Eulerian polynomial is $xA_{n}(x)$, for $n \ge 1$. We remark that our definition is the one that Euler defined originally in [6].

The first few terms of $A_{n}(x)$ are as follows:
\begin{align*}
&A_{1}(x)=1, A_{2}(x)=1+x, A_{3}(x)=1+4x+x^{2}, A_{4}(x)=1+11x+11x^{2}+x^{3}, \\
&A_{5}(x)=1+26x+66x^{2}+26x^{3}+x^{4},
A_{6}(x)=1+57x+302x^{2}+302x^{3}+57x^{4}+x^{5},\\
&A_{7}(x)=1+120x+1191x^{2}+2416x^{3}+1191x^{4}+120x^{5}+x^{6}.
\end{align*}
We list some of the properties of Eulerian polynomials (see [2,7,8,21]).\\
\begin{proposition}
\begin{align*}
&(a)\,\,A_{0}(x)=1,\,\, A_{n}(x)=\sum_{k=0}^{n-1}\binom{n}{k}A_{k}(x)(x-1)^{n-1-k},\,\, (n \ge 1), \\
&(b)\,\,\sum_{n=0}^{\infty}A_{n}(x)\frac{t^{n}}{n!}=\frac{1-x}{e^{t(x-1)}-x}, \\
&(c)\,\,\frac{A_{n}(x)}{(1-x)^{n+1}}=\sum_{j=0}^{\infty}x^{j}(j+1)^{n},\,\,(n \ge 0), \\
&(d)\,\,A_{0}(x)=1,\,\,A_{n}(x)=\big(1+(n-1)x\big)A_{n-1}(x)+x(1-x)A_{n-1}^{\prime}(x),\,\,(n \ge 1), \\
&(e)\,\,A_{n,k}=(k+1)A_{n-1,k}+(n-k)A_{n-1,k-1},\,\,(n \ge 2, 1 \le k \le n-1),\,\,A_{n,0}=1,\\
&(f)\,\,A_{n,k}=\sum_{i=0}^{k}(-1)^{i}(k+1-i)^{n}\binom{n+1}{i},\,\,(n \ge 0),\\
&(g) (\textrm{Worpitzky identity})\,\, x^{n}=\sum_{k=0}^{n}A_{n,k}\binom{x+k}{n}=\sum_{k=0}^{n}A_{n,k}\binom{x+n-k-1}{n},\,\,(n \ge 0), \\
&(h)\,\,\sum_{i=1}^{m} i^{n}x^{i}=\sum_{l=1}^{n}(-1)^{n+l}\binom{n}{l}\frac{x^{m+1}A_{n-l}(x)}{(x-1)^{n-l+1}}m^{l}+(-1)^{n}\frac{x(x^{m}-1)}{(x-1)^{n+1}}A_{n}(x),\,\,(m \ge 1, n \ge 0), \\
&(i)\,\,A_{n,k}=A_{n,n-1-k},\,\, (n \ge 1, 0 \le k \le n-1), \\
&(j)\,\,A_{n}(1)=\sum_{k=0}^{n}A_{n,k}=n!,\,\, (n \ge 0),\\
&(k) (\textrm{Frobenius})\,\,A_{n,k}=\sum_{l=1}^{k+1}l!S_{2}(n,l)\binom{n-l}{k-l+1}(-1)^{k-l-1},\,\,(n \ge 2, \, 1 \le k \le n-1),\\
&\quad A_{n,0}=1,\\
&(l)\,\,S_2(n,k)=\frac{1}{k!}\sum_{j=1}^{k}\binom{n-j}{k-j}A_{n,j-1}.
\end{align*}
\end{proposition}
Let $\bold{P}=\left\{p_{n}(x)\right\}_{n=0}^{\infty}$ be a sequence of polynomials such that deg\,$p_{n}(x)=n, p_{0}(x)=1$.
In [19],  Koutras defined,  motivated by the Worpitzky identity (g), the Eulerian numbers $A_{n,k}(\bold{P})$ associated with $\bold{P}$ by the coefficients in the expansion of $p_{n}(x)$ in terms of $\binom{x+n-k-1}{n}$. Thus 
\begin{align*}
p_{n}(x)=\sum_{k=0}^{n}A_{n,k}(\bold{P})\binom{x+n-k-1}{n}.
\end{align*}
He also defined the Eulerian polynomial associated with $\bold{P}$ by
\begin{align*} 
A_{n}(x;\bold{P})=\sum_{k=0}^{n}A_{n,k}(\bold{P})x^{k}.
\end{align*}
Some of the following are already stated in [19]. However, we state and prove them for the sake of completeness.
\begin{theorem}
Let $\bold{P}=\left\{p_{n}(x)\right\}_{n=0}^{\infty}$ be a sequence of polynomials such that deg\,$p_{n}(x)=n, p_{0}(x)=1$. Then we obtain\\
\begin{align*}
&(a)\,\,A_{n,k}(\bold{P})=\sum_{l=0}^{k}(-1)^{l}\binom{n+1}{l}p_{n}(k-l+1).\\
&(b)\,\,\frac{A_{n}(x;\bold{P})}{(1-x)^{n+1}}=\sum_{j=0}^{\infty}x^{j}p_{n}(j+1).\\
&(c)\,\, \textrm{Let $\bar{\bold{P}}=\left\{\bar{p}_{n}(x)\right\}$, where $\bar{p}_{n}(x)$ is the sequence of polynomials defined by}\\
&\quad \bar{p}_{0}(x)=1,\,\,\bar{p}_{n}(x)=xp_{n-1}(x),\,\,(n \ge 1). \,\,\textrm{Then we have}\\
&\quad A_{n+1}(x;\bar{\bold{P}})=(1+nx)A_{n}(x;\bold{P})+x(1-x)A_{n}^{\prime}(x;\bold{P}), \,\,(n \ge 0).\\
&(d)\,\, A_{n,k}(\bar{\bold{P}})=(k+1)A_{n-1,k}(\bold{P})+(n-k)A_{n-1,k-1}(\bold{P}), \,\,(n \ge 2, \,\,1 \le k \le n-1),\\
&\quad A_{n,0}(\bar{\bold{P}})=\bar{p}_{n}(1)=A_{n-1,0}(\bold{P})=p_{n-1}(1),\,\,(n \ge 1).\\
&(e)\,\, \textrm{Assume that $p_{n}(x)$ is Sheffer for the pair $(g(t),f(t))$. Then we have}\\
&\quad\sum_{n=0}^{\infty}A_{n}(x;\bold{P})\frac{t^{n}}{n!}=\frac{(1-x)e^{\bar{f}((1-x)t)}}{g(\bar{f}((1-x)t))}\frac{1}{1-xe^{\bar{f}((1-x)t)}}.\\
&(f)\,\, \textrm{Assume that $p_{n}(x)$ is Sheffer for the pair $(1,f(t))$, with $f(-t)=-f(t)$.  Then}\\
&\quad A_{n,k}(\bold{P})_=A_{n,n-1-k}(\bold{P}),\,\, (n \ge 1, 0 \le k \le n-1). \\
&\quad \textrm{For the assertions (g), (h), (i) and (j) below, assume that $p_{n}(0)=0$, for $ n \ge 1$.} \\
&\quad\textrm{Equivalently, $p_{n,0}=0$, for $n \ge 1.$  Then}\\
&(g)\,\,S_2(n,0;\bold{P})=0, \,\,\mathrm{deg}\,A_n(x;\bold{P}) \le n-1,  \,\,A_{n,n}(\bold{P})=0, \,\mathrm{for}\quad n \ge 1.\\
&\quad \textrm{More precisely, $A_n(x;\bold{P})=\sum_{k=1}^{n}(-1)^{n-k}p_{n,k}\, x^{n-1}+\cdots.$}\\
&(h)\,\,A_{n,k}(\bold{P})=\sum_{j=1}^{k+1}j!S_{2}(n,j;\bold{P})\binom{n-j}{k-j+1}(-1)^{k-j+1}.\\
&(i)\,\,S_2(n,k;\bold{P})=\frac{1}{k!}\sum_{j=1}^{k}\binom{n-j}{k-j}A_{n,j-1}(\bold{P}).\\
&(j)\,\,A_{n}(1;\bold{P})=p_{nn}n!. \\
&(k)\,\,p_{n}(x)=\sum_{k=0}^{n}A_{n,k}(\bold{P})\binom{x+n-k-1}{n}.
\end{align*}
\begin{proof}
(a) By using Vandermonde convolution formula (see [4,22,23]) 
\begin{align*}
\sum_{k=0}^{n}\binom{x}{k}\binom{y}{n-k}=\binom{x+y}{n},
\end{align*}
we have
\begin{align*}
&\sum_{l=0}^{k}(-1)^{l}\binom{n+1}{l}p_{n}(k-l+1)=\sum_{l=0}^{k}(-1)^{l}\binom{n+1}{l}\sum_{j=0}^{n}A_{n,j}(\bold{P})\binom{n+k-j-l}{n}\\
&=\sum_{l=0}^{k}(-1)^{l}\binom{n+1}{l}\sum_{j=0}^{k-l}A_{n,j}(\bold{P})\binom{n+k-j-l}{n}\\
&=\sum_{j=0}^{k}A_{n,j}(\bold{P})(-1)^{k-j}\sum_{l=0}^{k-j}(-1)^{k-j-l}\binom{n+1}{l}\binom{n+k-j-l}{n}\\
&=\sum_{j=0}^{k}A_{n,j}(\bold{P})(-1)^{k-j}\sum_{l=0}^{k-j}\binom{n+1}{l}\binom{-n-1}{k-j-l}\\
&=\sum_{j=0}^{k}A_{n,j}(\bold{P})(-1)^{k-j}\delta_{k,j}=A_{n,k}(\bold{P}).
\end{align*}
(b)\,\,By making use of (a), we see that
\begin{align*}
A_{n}(x;\bold{P})&=\sum_{k=0}^{n}A_{n,k}(\bold{P})x^k=
\sum_{k=0}^{n}\sum_{l=0}^{k}(-1)^{l}\binom{n+1}{l}p_{n}(k-l+1)x^{k}\\
&=\sum_{l=0}^{\infty}(-1)^{l}\binom{n+1}{l}x^{l}\sum_{j=0}^{\infty}p_{n}(j+1)x^{j}\\
&=(1-x)^{n+1}\sum_{j=0}^{\infty}p_{n}(j+1)x^{j}.
\end{align*}
(c)\,\, This follows by multiplying the identity in (b) by $x$ and then differentiating the resulting identity.\par
(d)\,\,From (c), it is easy to see that we have\\
\begin{align*}
\sum_{k=0}^{n+1}A_{n+1,k}(\bar{\bold{P}})x^{k}&=\sum_{k=0}^{n}A_{n,k}(\bold{P})x^{k}+\sum_{k=1}^{n}k A_{n,k}(\bold{P})x^{k}\\
&+\sum_{k=1}^{n+1}n A_{n,k-1}(\bold{P})x^{k}-\sum_{k=2}^{n+1}(k-1)A_{n,k-1}(\bold{P})x^{k}\\
&=\sum_{k=0}^{n}A_{n,k}(\bold{P})x^{k}+\sum_{k=1}^{n}k A_{n,k}(\bold{P})x^{k}\\
&+\sum_{k=1}^{n}n A_{n,k-1}(\bold{P})x^{k}-\sum_{k=1}^{n}(k-1)A_{n,k-1}(\bold{P})x^{k}.
\end{align*}
(e)\,\,By making use of (b), we have
\begin{align*}
\sum_{n=0}^{\infty}A_{n}(x;\bold{P})\frac{t^{n}}{n!}&=\sum_{n=0}^{\infty}(1-x)^{n+1}\sum_{j=0}^{\infty}p_{n}(j+1)x^{j}\frac{t^{n}}{n!}\\
&=(1-x)\sum_{j=0}^{\infty}x^{j}\sum_{n=0}^{\infty}p_{n}(j+1)\frac{((1-x)t)^{n}}{n!}\\
&=(1-x)\sum_{j=0}^{\infty}x^{j}\frac{1}{g(\bar{f}((1-x)t))}e^{(j+1)\bar{f}((1-x)t)},
\end{align*}
which gives the desired result.\\
(f)\,\,Note first that $\bar{f}(-t)=-\bar{f}(t).$ Here the assertion is equivalent to showing that $A_{n}(x;\bold{P})=x^{n-1}A_{n}(\frac{1}{x};\bold{P})$, which means it is a Euler-Frobenius polynomial. By using (c), we see that
\begin{align*}
\sum_{n=1}^{\infty}A_{n}(x;\bold{P})\frac{t^{n}}{n!}=\frac{(1-x)e^{\bar{f}((1-x)t)}}{1-xe^{\bar{f}((1-x)t)}}-1=\frac{e^{\bar{f}((1-x)t)}-1}{1-xe^{\bar{f}((1-x)t)}}.
\end{align*}
By invoking $\bar{f}(-t)=-\bar{f}(t)$, we get the desired result, as we see from
\begin{align*}
\sum_{n=1}^{\infty}x^{n-1}A_{n}(\frac{1}{x};\bold{P})\frac{t^{n}}{n!}=\frac{e^{\bar{f}((x-1)t)}-1}{x-e^{\bar{f}((x-1)t)}}=\frac{e^{\bar{f}((1-x)t)}-1}{1-xe^{\bar{f}((1-x)t)}}.
\end{align*}
(g)\,\, Let $n \ge 1$. By invoking (b) and Theorem 3.2 (c), we see that
\begin{align*}
A_{n}(x;\bold{P})&=(1-x)^{n+1}\sum_{j=0}^{\infty}x^{j}\sum_{k=1}^{n}p_{nk}(j+1)^{k}\\
&=\sum_{k=1}^{n}p_{nk}(1-x)^{n-k}(1-x)^{k+1}\sum_{j=0}^{\infty}x^{j}(j+1)^{k}\\
&=\sum_{k=1}^{n}p_{nk}(1-x)^{n-k}A_{k}(x)\\
&=\sum_{k=1}^{n}p_{nk}(-1)^{n-k}x^{n-1}+\cdots,
\end{align*}
where we used the fact that $A_{k}(x)$ is a monic polynomial of degree $k-1$, for $k \ge 1$.\\
This shows in particular that $A_{n,n}(\bold{P})=0$. The fact that $S_{2}(n,0;\bold{P})=0$ follows from the definition \eqref{5C}. \\
(h)\,\,To prove this, we may assume that $n \ge 1$. By using (b) and \eqref{5C}, we have
\begin{align}
A_{n}(x;\bold{P})&=(1-x)^{n+1}\sum_{k=0}^{\infty}x^{k}p_{n}(k+1)\label{6C}\\
&=(1-x)^{n+1}\sum_{j=0}^{n}S_{2}(n,j;\bold{P})\sum_{k=0}^{\infty}(k+1)_{j}x^{k},\nonumber
\end{align}
where we note that
\begin{align}
\sum_{k=0}^{\infty}(k+1)_{j}x^{k}&=\sum_{k=j-1}^{\infty}(k+1)_{j}x^{k}=\sum_{k=0}^{\infty}(k+j)_{j}x^{k+j-1}
=j!x^{j-1}(1-x)^{-(j+1)}.\label{7C}
\end{align}
Thus we see from \eqref{6C} and \eqref{7C} that
\begin{align*}
A_{n}(x;\bold{P})&=\sum_{j=0}^{n}j!x^{j-1}S_{2}(n,j;\bold{P})(1-x)^{n-j}\\
&=\sum_{j=0}^{n}j!x^{j-1}S_{2}(n,j;\bold{P})\sum_{m=0}^{\infty}\binom{n-j}{m}(-1)^{m}x^{m}\\
&=\sum_{j=0}^{n}\sum_{m=0}^{\infty}j!S_{2}(n,j;\bold{P})\binom{n-j}{m}(-1)^{m}x^{j+m-1}\\
&=\sum_{k=0}^{n}\sum_{j=0}^{k}j!S_{2}(n,j;\bold{P})\binom{n-j}{k-j}(-1)^{k-j}x^{k-1}\\
&=\sum_{k=1}^{n}\sum_{j=1}^{k}j!S_{2}(n,j;\bold{P})\binom{n-j}{k-j}(-1)^{k-j}x^{k-1}\\
&=\sum_{k=0}^{n-1}\sum_{j=1}^{k+1}j!S_{2}(n,j;\bold{P})\binom{n-j}{k-j+1}(-1)^{k-j+1}x^{k},
\end{align*}
from which the result follows.\\
(i)\,\,By using (h), the right hand side of (i) is 
\begin{align*}
\frac{1}{k!}\sum_{j=1}^{k}\binom{n-j}{k-j}A_{n,j-1}(\bold{P})&=\frac{1}{k!}\sum_{j=1}^{k}\binom{n-j}{k-j}\sum_{l=1}^{j}l!S_{2}(n,l;\bold{P})\binom{n-l}{j-l}(-1)^{j-l}\\
&=\frac{1}{k!}\sum_{l=1}^{k}l!S_{2}(n,l;\bold{P})\sum_{j=l}^{k}\binom{n-j}{k-j}\binom{n-l}{j-l}(-1)^{j-l}\\
&=\frac{1}{k!}\sum_{l=1}^{k}l!S_{2}(n,l;\bold{P})\sum_{j=l}^{k}\binom{k-l}{k-j}\binom{n-l}{n-k}(-1)^{j-l}\\
&=\frac{1}{k!}\sum_{l=1}^{k}l!S_{2}(n,l;\bold{P})\binom{n-l}{n-k}\sum_{j=l}^{k}\binom{k-l}{k-j}(-1)^{j-l}\\
&=\frac{1}{k!}\sum_{l=1}^{k}l!S_{2}(n,l;\bold{P})\binom{n-l}{n-k}(-1)^{k-l}\sum_{j=0}^{k-l}\binom{k-l}{j}(-1)^{j}\\
&=\frac{1}{k!}\sum_{l=1}^{k}l!S_{2}(n,l;\bold{P})\binom{n-l}{n-k}(-1)^{k-l}\delta_{k,l}\\
&=S_{2}(n,k;\bold{P}).
\end{align*}
(j)\,\, As we saw in the proof of (g) and by using Proposition 3.2 (j), we have
\begin{align*}
A_{n}(x;\bold{P})=\sum_{k=1}^{n}p_{nk}(1-x)^{n-k}A_{k}(x),
\end{align*}
from which the assertion follows. \\
(k)\,\,This is just the definition in \eqref{5C}.
\end{proof}
\end{theorem}

\section{Examples on Stirling numbers of the second kind associated with sequences of polynomials}
(a) Let $\bold{P}=\left\{x^{n}\right\}$. Then $x^{n} \sim (1,t)$.
By the definition in \eqref{5C} and Theorem 3.1, we have (see [4,22,23])
\begin{align*}
S_{2}(n,k;\bold{P})=S_{2}(n,k), \quad \sum_{n=k}^{\infty}S_{2}(n,k)\frac{t^n}{n!}=\frac{1}{k!}(e^{t}-1)^{k}.
\end{align*}
(b) Let $\bold{P}=\left\{(x)_{n,\lambda}\right\}$ be the sequence of generalized falling factorials. Then $(x)_{n,\lambda} \sim \big(1, \frac{1}{\lambda}(e^{\lambda t}-1)\big)$. By \eqref{4C} and \eqref{5C}, we have (see [9,10,12,14,16-18])
\begin{align*}
S_{2}(n,k;\bold{P})=S_{2,\lambda}(n,k).
\end{align*}
Further, by Theorem 3.1, we get
\begin{align*}
\sum_{n=k}^{\infty}S_{2,\lambda}(n,k)\frac{t^{n}}{n!}=\frac{1}{k!}\big((1+\lambda t)^{\frac{1}{\lambda}}-1\big)^{k}=\frac{1}{k!}(e_{\lambda}(t)-1)^{k},
\end{align*}
where we use the standard notations for degenerate exponentials:
\begin{align*}
e_{\lambda}^{x}(t)=(1+\lambda t)^{\frac{x}{\lambda}}, \quad e_{\lambda}(t)=e_{\lambda}^{1}(t)=(1+\lambda t)^{\frac{1}{\lambda}}.
\end{align*}
(c) The rising factorials $\langle x \rangle_{n}$ are defined by 
\begin{align*}
\langle x \rangle_{n}=x(x+1)\cdots(x+n-1), \,\,(n \ge 1), \quad \langle x \rangle_{0}=1.
\end{align*}
Let $\bold{P}=\left\{\langle{x \rangle}_{n}\right\}$. Then $\langle{x \rangle}_{n} \sim (1, 1-e^{-t}),$ and $\langle{x \rangle}_{n}=\sum_{k=0}^{n}L(n,k)(x)_{k}$. Here $L(n,k)=\binom{n-1}{k-1}\frac{n!}{k!}$ are the (unsigned) Lah numbers and (see [4,11,22,23])
\begin{align*}
S_{2}(n,k;\bold{P})=L(n,k).
\end{align*}
An expression of the Lah numbers can be derived from Theorem 3.1. Indeed, we have
\begin{align*}
\sum_{n=k}^{\infty}L(n,k)\frac{t^{n}}{n!}&=\frac{1}{k!}\big(\frac{t}{1-t}\big)^{k}=\frac{1}{k!}\big(e^{-\log (1-t)}-1\big)^{k}\\
&=\sum_{l=k}^{\infty}S_{2}(l,k)(-1)^{l}\frac{1}{l!}\big(\log (1-t) \big)^{l}\\
&=\sum_{n=k}^{\infty}\sum_{l=k}^{n}(-1)^{n-l}S_{1}(n,l)S_{2}(l,k)\frac{t^{n}}{n!}.
\end{align*}
Hence we have
\begin{align*}
L(n,k)=\sum_{l=k}^{n}(-1)^{n-l}S_{1}(n,l)S_{2}(l,k).
\end{align*}
(d) The generalized rising factorials are defined by
\begin{align*}
\langle x \rangle_{n,\lambda}=x(x+\lambda)\cdots(x+(n-1)\lambda), \,\,(n \ge 1), \quad \langle x \rangle_{0,\lambda}=1.
\end{align*}
Let $\bold{P}=\left\{\langle{x \rangle}_{n,\lambda}\right\}$. Then $\langle{x \rangle}_{n,\lambda} \sim \big(1, \frac{1}{\lambda}(1-e^{-\lambda t})\big)$. The degenerate Lah numbers are defined by $\langle{x \rangle}_{n,\lambda}=\sum_{k=0}^{n}L_{\lambda}(n,k)(x)_{k}$, and hence (see [14])
\begin{align*}
S_{2}(n,k;\bold{P})=L_{\lambda}(n,k).
\end{align*}
An expression of the degenerate Lah numbers can be derived from Theorem 3.1. Indeed, we have
\begin{align*}
\sum_{n=k}^{\infty}L_{\lambda}(n,k)\frac{t^{n}}{n!}&=\frac{1}{k!}\big(e^{-\frac{1}{\lambda}\log (1-\lambda t)}-1\big)^{k}\\
&=\sum_{l=k}^{\infty}S_{2}(l,k)(-\frac{1}{\lambda})^{l}\frac{1}{l!}\big(\log (1-\lambda t) \big)^{l}\\
&=\sum_{n=k}^{\infty}\sum_{l=k}^{n}(-\lambda)^{n-l}S_{1}(n,l)S_{2}(l,k)\frac{t^{n}}{n!}.
\end{align*}
Hence we have
\begin{align*}
L_{\lambda}(n,k)=\sum_{l=k}^{n}(-\lambda)^{n-l}S_{1}(n,l)S_{2}(l,k).
\end{align*}
(e) The central factorials $x^{[n]}$ are defined by (see [3,4,15,23])
\begin{align*}
x^{[n]}=x(x+\frac{1}{2}n-1)_{n-1}, \,\,(n \ge 1),\quad x^{[0]}=1,
\end{align*}
where $x^{[n]} \sim (1, f(t)=e^{\frac{t}{2}}-e^{-\frac{t}{2}})$. Here we note that
\begin{align*}
\bar{f}(t)=2\log\bigg(\frac{t+\sqrt{t^{2}+4}}{2}\bigg)=\log \big(1+\frac{t}{2}(t+\sqrt{t^{2}+4})\big).
\end{align*}
Carlitz-Riordan and Riordan (see [3,4,15,23]) discussed the numbers $T_{1}(n,k)$ and $T_{2}(n,k)$ which are respectively called the central factorial numbers of the first kind and of the second kind, and defined by
\begin{align*}
&\frac{1}{k!}\bigg(2\log\bigg(\frac{t+\sqrt{t^{2}+4}}{2}\bigg)\bigg)^{k}=\frac{1}{k!}\Big(\log \big(1+\frac{t}{2}(t+\sqrt{t^{2}+4}\big)\big)\Big)^{k} \\
&=\sum_{n=k}^{\infty}T_{1}(n,k)\frac{t^{n}}{n!},\quad x^{[n]}=\sum_{k=0}^{n}T_{1}(n,k)x^{k},\\
&\frac{1}{k!}\big(e^{\frac{t}{2}}-e^{-\frac{t}{2}}\big)^{k}=\sum_{n=k}^{\infty}T_{2}(n,k)\frac{t^{n}}{n!},\quad x^{n}=\sum_{k=0}^{n}T_{2}(n,k)x^{[k]}. \\
\end{align*}
Let $\bold{P}=\left\{x^{[n]}\right\}$. From Theorem 3.1, we get
\begin{align*}
S_{2}(n,k;\bold{P})=\sum_{l=k}^{n}S_{2}(l,k)T_{1}(n,l),\quad \sum_{n=k}^{\infty}S_{2}(n,k;\bold{P})\frac{t^{n}}{n!}=\frac{1}{k!}\big(\frac{t}{2}(t+\sqrt{t^{2}+4})\big)^{k}.
\end{align*}
(f) The central Bell polynomials $\mathrm{Bel}_{n}^{(c)}(x)$ and the central Bell numbers $\mathrm{Bel}_{n}^{(c)}$ are respectively defined by (see [15])
\begin{align*}
e^{x(e^{\frac{t}{2}}-e^{-\frac{t}{2}})}=\sum_{n=0}^{\infty}\mathrm{Bel}_{n}^{(c)}(x)\frac{t^{n}}{n!}, \quad \mathrm{Bel}_{n}^{(c)}=\mathrm{Bel}_{n}^{(c)}(1).
\end{align*}
Then we see that 
\begin{align*}
\mathrm{Bel}_{n}^{(c)}(x) \sim \Big(1, f(t)=2\log\bigg(\frac{t+\sqrt{t^{2}+4}}{2}\bigg)\Big), \quad \mathrm{Bel}_{n}^{(c)}(x)=\sum_{k=0}^{n}T_{2}(n,k) x^{k}.
\end{align*}
Let $\bold{P}=\left\{\mathrm{Bel}_{n}^{(c)}(x)\right\}$. Then, by Theorem 3.1 and (e), we see that
\begin{align*}
S_{2}(n,k;\bold{P})=\sum_{l=k}^{n}S_{2}(l,k)T_{2}(n,l),\quad \sum_{n=k}^{\infty}S_{2}(n,k;\bold{P})\frac{t^{n}}{n!}=\frac{1}{k!}\big(e^{(e^{\frac{t}{2}}-e^{-\frac{t}{2}})}-1\big)^{k}.
\end{align*}
(g) The numbers $T_{1,\lambda}(n,k)$ and $T_{2,\lambda}(n,k)$ are respectively called the degenerate central factorial numbers of the first kind and of the second kind, and defined by (see [16])
\begin{align*}
&\frac{1}{k!}\left\{\frac{1}{\lambda}\bigg(\bigg(\frac{t+\sqrt{t^{2}+4}}{2}\bigg)^{2 \lambda}-1\bigg)\right\}^{k}=\frac{1}{k!}\bigg(\log_{\lambda}\bigg(\frac{t+\sqrt{t^{2}+4}}{2}\bigg)^{2}\bigg)^{k}=\sum_{n=k}^{\infty}T_{1,\lambda}(n,k)\frac{t^{n}}{n!},\\
&x^{[n]}=\sum_{k=0}^{n}T_{1,\lambda}(n,k)(x)_{k, \lambda},\\
&\frac{1}{k!}\big(e_{\lambda}^{\frac{1}{2}}(t)-e_{\lambda}^{-\frac{1}{2}}(t)\big)^{k}=\sum_{n=k}^{\infty}T_{2,\lambda}(n,k)\frac{t^{n}}{n!},\quad (x)_{n,\lambda}=\sum_{k=0}^{n}T_{2,\lambda}(n,k)x^{[k]}.
\end{align*}
The degenerate central Bell polynomials $\mathrm{Bel}_{n,\lambda}^{(c)}(x)$ and the degenerate central Bell numbers $\mathrm{Bel}_{n,\lambda}^{(c)}$ are respectively defined by (see [16])
\begin{align*}
e^{x(e_{\lambda}^{\frac{1}{2}}(t)-e_{\lambda}^{-\frac{1}{2}}(t))}=\sum_{n=0}^{\infty}\mathrm{Bel}_{n,\lambda}^{(c)}(x)\frac{t^{n}}{n!}, \quad \mathrm{Bel}_{n,\lambda}^{(c)}=\mathrm{Bel}_{n,\lambda}^{(c)}(1).
\end{align*}
Then we see that 
\begin{align*}
\mathrm{Bel}_{n,\lambda}^{(c)}(x) \sim \Big(1, f(t)=\frac{1}{\lambda}\bigg(\bigg(\frac{t+\sqrt{t^{2}+4}}{2}\bigg)^{2 \lambda}-1\bigg)\Big),\quad
\mathrm{Bel}_{n,\lambda}^{(c)}(x)=\sum_{k=0}^{n}T_{2,\lambda}(n,k) x^{k}.
\end{align*}
Let $\bold{P}=\left\{\mathrm{Bel}_{n}^{(c)}(x)\right\}$. Then, by Theorem 3.1 and (e), we see that
\begin{align*}
S_{2}(n,k;\bold{P})=\sum_{l=k}^{n}S_{2}(l,k)T_{2,\lambda}(n,l),\quad \sum_{n=k}^{\infty}S_{2}(n,k;\bold{P})\frac{t^{n}}{n!}=\frac{1}{k!}\big(e^{(e_{\lambda}^{\frac{1}{2}}(t)-e_{\lambda}^{-\frac{1}{2}}(t))}-1\big)^{k}.
\end{align*}
(h) Let $x^{[n,\lambda]}$ be defined by
\begin{align*}
x^{[n,\lambda]}=x\big(x+(\frac{1}{2}n-1)\lambda\big)_{n-1,\lambda}, \,\,(n \ge 1),\quad x^{[0,\lambda]}=1.
\end{align*}
where $x^{[n,\lambda]} \sim \big(1, f(t)=\frac{1}{\lambda}(e^{\frac{\lambda t}{2}}-e^{-\frac{\lambda t}{2}})\big)$. Here we note that
\begin{align*}
\bar{f}(t)=\frac{2}{\lambda}\log\bigg(\frac{\lambda t+\sqrt{\lambda^{2}t^{2}+4}}{2}\bigg)=\frac{1}{\lambda}\log \big(1+\frac{\lambda t}{2}(\lambda t+\sqrt{\lambda^{2}t^{2}+4})\big).
\end{align*}
As analogies to the numbers $T_{1,\lambda}(n,k)$ and $T_{2,\lambda}(n,k)$ mentioned in (e) above, we may introduce the numbers $R_{1,\lambda}(n,k)$ and $R_{2,\lambda}(n,k)$ which are defined by
\begin{align*}
&x^{[n,\lambda]}=\sum_{k=0}^{n}R_{1,\lambda}(n,k)x^{k},\quad \sum_{n=k}^{\infty}R_{1,\lambda}(n,k)\frac{t^{n}}{n!}=\frac{1}{k!}\bigg(\frac{2}{\lambda}\log\bigg(\frac{\lambda t+\sqrt{\lambda^{2}t^{2}+4}}{2}\bigg)\bigg)^{k},\\
&x^{n}=\sum_{k=0}^{n}R_{2,\lambda}(n,k)x^{[k,\lambda]}, \quad \sum_{n=k}^{\infty}R_{2,\lambda}(n,k)\frac{t^{n}}{n!}=\frac{1}{k!}\Big(\frac{1}{\lambda}\big(e^{\frac{\lambda t}{2}}-e^{-\frac{\lambda t}{2}}\big)\Big)^{k}.
\end{align*}
Let $\bold{P}=\left\{x^{[n,\lambda]}\right\}$. From Theorem 3.1, we get
\begin{align*}
S_{2}(n,k;\bold{P})=\sum_{l=k}^{n}S_{2}(l,k)R_{1,\lambda}(n,l),\,\, \sum_{n=k}^{\infty}S_{2}(n,k;\bold{P})\frac{t^{n}}{n!}=\frac{1}{k!}\bigg(\bigg(\frac{\lambda t+\sqrt{\lambda^{2}t^{2}+4}}{2}\bigg)^{\frac{2}{\lambda}}-1\bigg)^{k}.
\end{align*}
(i) Let $\bold{P}=\left\{B_{n}^{L}(x)\right\}$ be the sequence of Lah-Bell polynomials. Then it is given by $e^{x(\frac{t}{1-t})}=\sum_{n=0}^{\infty}B_{n}^{L}(x)\frac{t^n}{n!}$, so that $B_{n}^{L}(x) \sim (1,\frac{t}{1+t})$, and $B_{n}^{L}(x)=\sum_{k=0}^{n}L(n,k)x^{k}$ (see [11]).  By Theorem 3.1 and (c), we have
\begin{align*}
\sum_{n=k}^{\infty}S_{2}(n,k;\bold{P})\frac{t^{n}}{n!}&=\frac{1}{k!}(e^{\frac{t}{1-t}}-1)^{k}=\sum_{l=k}^{\infty}S_{2}(l,k)\frac{1}{l!}\big(\frac{t}{1-t}\big)^{l}\\
&=\sum_{l=k}^{\infty}S_{2}(l,k)\sum_{n=l}^{\infty}L(n,l)\frac{t^{n}}{n!}\\
&=\sum_{n=k}^{\infty}\sum_{l=k}^{n}L(n,l)S_{2}(l,k)\frac{t^{n}}{n!}.
\end{align*}
Thus we obtain $S_{2}(n,k;\bold{P})=\sum_{l=k}^{n}L(n,l)S_{2}(l,k)$. \\
(j) Let $\bold{P}=\left\{B_{n,\lambda}^{L}(x)\right\}$ be the sequence of degenerate Lah-Bell polynomials. Then it is given by $e_{\lambda}^{x}(\frac{t}{1-t})=\sum_{n=0}^{\infty}B_{n,\lambda}^{L}(x)\frac{t^n}{n!}$. Thus we see that $B_{n,\lambda}^{L}(x) \sim (1,\frac{e^{\lambda t}-1}{\lambda+e^{\lambda t}-1})$, and $B_{n,\lambda}^{L}(x)=\sum_{k=0}^{n}L(n,k)(x)_{k,\lambda}$ (see [14]).  By Theorem 3.1 and recalling from (d) that $\sum_{n=k}^{\infty}L_{\lambda}(n,k)\frac{t^{n}}{n!}=\frac{1}{k!}\big((1-\lambda t)^{-\frac{1}{\lambda}}-1\big)^{k}$, we have
\begin{align*}
\sum_{n=k}^{\infty}S_{2}(n,k;\bold{P})\frac{t^{n}}{n!}&=\frac{1}{k!}(e^{\frac{1}{\lambda}\log(1+\frac{\lambda t}{1-t})}-1)^{k}\\
&=\sum_{m=k}^{\infty}L_{-\lambda}(m,k)\frac{1}{m!}\big(\frac{t}{1-t}\big)^{m}\\
&=\sum_{m=k}^{\infty}L_{-\lambda}(m,k)\sum_{n=m}^{\infty}L(n,m)\frac{t^{n}}{n!}\\
&=\sum_{n=k}^{\infty}\sum_{m=k}^{n}L(n,m)L_{-\lambda}(m,k)\frac{t^{n}}{n!}
\end{align*}
Thus we obtain $S_{2}(n,k;\bold{P})=\sum_{m=k}^{n}L(n,m)L_{-\lambda}(m,k).$ \\
(k) Let $\bold{P}=\left\{\mathrm{Bel}_{n}(x)\right\}$ be the sequence of Bell polynomials given by $e^{x(e^{t}-1)}=\sum_{n=0}^{\infty}\mathrm{Bel}_{n}(x)\frac{t^{n}}{n!}$. Then $\mathrm{Bel}_{n}(x) \sim (1, \log (1+t))$, and $\mathrm{Bel}_{n}(x)=\sum_{k=0}^{\infty}S_{2}(n,k)x^{k}$ (see [4,22,23]).  Now, by Theorem 3.1, we have
\begin{align*}
\sum_{n=k}^{\infty}S_{2}(n,k;\bold{P})\frac{t^{n}}{n!}&=\frac{1}{k!}\big(e^{e^t-1}-1\big)^{k}=\sum_{l=k}^{\infty}S_{2}(l,k)\frac{1}{l!}(e^{t}-1)^{l}\\
&=\sum_{n=k}^{\infty}\sum_{l=k}^{n}S_{2}(l,k)S_2(n,l)\frac{t^{n}}{n!}.
\end{align*}
Thus we obtain
\begin{align*}
S_{2}(n,k;\bold{P})=\sum_{l=k}^{n}S_2(n,l)S_{2}(l,k).
\end{align*}
(l) Let $\bold{P}=\left\{\mathrm{Bel}_{n,\lambda}(x)\right\}$ be the sequence of partially degenerate Bell polynomials given by  $e^{x(e_{\lambda}(t)-1)}=\sum_{n=0}^{\infty}\mathrm{Bel}_{n,\lambda}(x)\frac{t^{n}}{n!}$. Then $\mathrm{Bel}_{n,\lambda}(x) \sim (1, \log_{\lambda}(1+t))$, and $\mathrm{Bel}_{n,\lambda}(x)=\sum_{k=0}^{n}S_{2,\lambda}(n,k)x^{k}$ (see [17]). Here $\log_{\lambda}(t)=\frac{1}{\lambda}(t^{\lambda}-1)$ is the degenerate logarithm which is the compositional inverse to the degenerate exponential $e_{\lambda}(t)$. From Theorem 3.1 and proceeding just in the case of (k), it is immediate to see that
\begin{align*}
S_{2}(n,k;\bold{P})=\sum_{l=k}^{n}S_{2,\lambda}(n,l)S_{2}(l,k).
\end{align*}
There are several other degenerate versions of Bell polynomials, one of which is the so called fully degenerate Bell polynomials $\phi_{n,\lambda}(x)$ given by $e_{\lambda}^{x}(e_{\lambda}(t)-1)=\sum_{n=0}^{\infty}\phi_{n,\lambda}(x)\frac{t^{n}}{n!}.$ Then $\phi_{n,\lambda}(x) \sim \big(1, \log_{\lambda}(1+\frac{1}{\lambda}(e^{\lambda t}-1))\big)$, and $\phi_{n,\lambda}(x)=\sum_{k=0}^{n}S_{2,\lambda}(n,k)(x)_{k,\lambda}$ (see [18]).  Now,  from Theorem 3.1 and recalling from (d) that $\sum_{n=k}^{\infty}L_{\lambda}(n,k)\frac{t^{n}}{n!}=\frac{1}{k!}\big((1-\lambda t)^{-\frac{1}{\lambda}}-1\big)^{k}$, we get
\begin{align*}
\sum_{n=k}^{\infty}S_{2}(n,k;\bold{P})\frac{t^{n}}{n!}&=\frac{1}{k!}\big(e^{\frac{1}{\lambda}\log(1+\lambda(e_{\lambda}(t)-1))}-1\big)^{k}\\
&=\frac{1}{k!}\big((1+\lambda(e_{\lambda}(t)-1))^{\frac{1}{\lambda}}-1\big)^{k}\\
&=\sum_{m=k}^{\infty}L_{-\lambda}(m,k)\frac{1}{m!}(e_{\lambda}(t)-1)^{m}\\
&=\sum_{n=k}^{\infty}\sum_{m=k}^{n}S_{2,\lambda}(n,m)L_{-\lambda}(m,k)\frac{t^{n}}{n!},
\end{align*}
which verifies that we have $S_{2}(n,k;\bold{P})=\sum_{m=k}^{n}S_{2,\lambda}(n,m)L_{-\lambda}(m,k).$\\
(m) Let $\bold{P}=\left\{M_{n}(x)\right\}$ be the sequence of Mittag-Leffler polynomials. That is, $M_{n}(x) \sim (1,f(t)=\frac{e^{t}-1}{e^{t}+1})$, with $\bar{f}(t)=\log \big(\frac{1+t}{1-t}\big)$ (see [23]).
By Theorem 3.1, we have
\begin{align*}
\sum_{n=k}^{\infty}S_{2}(n,k;\bold{P})\frac{t^{n}}{n!}&=\frac{1}{k!}\big(e^{\log(\frac{1+t}{1-t})}-1\big)^{k}=2^{k}\frac{1}{k!}\big(\frac{t}{1-t}\big)^{k}\\
&=\sum_{n=k}^{\infty}2^{k}L(n,k)\frac{t^{n}}{n!}.
\end{align*}
Thus we have
\begin{align*}
S_{2}(n,k;\bold{P})=2^{k}L(n,k).
\end{align*}
(n) Let $\bold{P}=\left\{L_{n}(x)\right\}$ be the sequence of Laguerre polynomials of order -1 (see [23]).
That is, $L_{n}(x) \sim (1, f(t)=\frac{t}{t-1})$, with $\bar{f}(t)=\frac{t}{t-1}$. By Theorem 3.1, we have
\begin{align*}
\sum_{n=k}^{\infty}S_{2}(n,k;\bold{P})\frac{t^{n}}{n!}&=\frac{1}{k!}\big(e^{\frac{t}{t-1}}-1\big)^{k}=\sum_{l=k}^{\infty}(-1)^{l}S_{2}(l,k)\frac{1}{l!}\big(\frac{t}{1-t}\big)^{l}\\
&=\sum_{l=k}^{\infty}(-1)^{l}S_{2}(l,k)\sum_{n=l}^{\infty}L(n,k)\frac{t^{n}}{n!}\\
&=\sum_{n=k}^{\infty}\sum_{l=k}^{n}(-1)^{l}L(n,l)S_{2}(l,k)\frac{t^{n}}{n!}.
\end{align*}
Thus we have $S_{2}(n,k;\bold{P})=\sum_{l=k}^{n}(-1)^{l}L(n,l)S_{2}(l,k)$. \\
(o) Let $\bold{P}=\left\{B_{n}(x)\right\}$ be the sequence of Bernoulli polynomials. Then $B_{n}(x) \sim \big(\frac{e^{t}-1}{t},t\big)$, and $B_{n}(x)=\sum_{k=0}^{n}\binom{n}{k}B_{n-k}x^{k}$ (see [23]).  By Theorem 3.1, we have
\begin{align*}
S_{2}(n,k;\bold{P})&=\sum_{l=k}^{n}S_{2}(l,k)\binom{n}{l}B_{n-l}=\sum_{l=0}^{n-k}\binom{n}{l}S_{2}(n-l,k)B_{l},\\
&\sum_{n=k}^{\infty}S_{2}(n,k;\bold{P})\frac{t^{n}}{n!}=\frac{t}{e^{t}-1}\frac{1}{k!}(e^{t}-1)^{k}.
\end{align*}
(p) Let $\bold{P}=\left\{E_{n}(x)\right\}$ be the sequence of Euler polynomials. Then $E_{n}(x) \sim \big(\frac{e^{t}+1}{2},t\big)$ (see [23]).
Analogously to (c), we see that
\begin{align*}
&S_{2}(n,k;\bold{P})=\sum_{l=0}^{n-k}\binom{n}{l}S_{2}(n-l,k)E_{l},\\
&\sum_{n=k}^{\infty}S_{2}(n,k;\bold{P})\frac{t^{n}}{n!}=\frac{2}{e^{t}+1}\frac{1}{k!}(e^{t}-1)^{k}.
\end{align*}
(q) Let $\bold{P}=\left\{(rx+s)_{n}\right\}$, with $r \ne 0$. Then
\begin{align*}
(rx+s)_{n}=\sum_{k=0}^{n}G(n,k;r,s)(x)_{k}, 
\end{align*}
where $G(n,k;r,s)$ are called the Gould-Hopper numbers (see [19]).  Thus
\begin{align*}
S_{2}(n,k;\bold{P})=G(n,k;r,s).
\end{align*}
Noting that $(rx+s)_{n} \sim \big(e^{-\frac{s}{r} t}, e^{\frac{t}{r}}-1\big)$, we can derive an expression of Gould-Hopper numbers from Theorem 3.1. Indeed, we get
\begin{align*}
&\sum_{n=k}^{\infty}G(n,k;r,s)\frac{t^{n}}{n!}=(1+t)^{s}\frac{1}{k!}\big(e^{r \log(1+t)}-1\big)^{k}\\
&=\sum_{j=0}^{\infty}(s)_{j}\frac{t^{j}}{j!}\sum_{l=k}^{\infty}\bigg(\sum_{m=k}^{l}r^{m}S_1(l,m)S_2(m,k)\bigg)\frac{t^{l}}{l!}\\
&=\sum_{n=k}^{\infty}\sum_{l=k}^{n}\sum_{m=k}^{l}\binom{n}{l}r^{m}(s)_{n-l}S_{1}(l,m)S_{2}(m,k) \frac{t^{n}}{n!}.
\end{align*}
Thus we have
\begin{align*}
G(n,k;r,s)=\sum_{l=k}^{n}\sum_{m=k}^{l}\binom{n}{l}r^{m}(s)_{n-l}S_{1}(l,m)S_{2}(m,k).
\end{align*}
(r) Let $\bold{P}=\left\{b_{n}(x)\right\}$ be the sequence of Bernoulli polynomials of the second kind. Then $b_{n}(x) \sim \big(\frac{t}{e^{t}-1}, e^{t}-1\big)$, so that $\sum_{n=0}^{\infty}b_{n}(x)\frac{t^{n}}{n!}=\frac{t}{\log (1+t)}(1+t)^{x}$, (see [23]).  From Theorem 3.1, we get
\begin{align*}
\sum_{n=k}^{\infty}S_{2}(n,k;\bold{P})\frac{t^{n}}{n!}&=\frac{t}{\log (1+t)}\frac{1}{k!}\big(e^{\log (1+t)}-1\big)^{k}\\
&=\frac{1}{k!}t^k\sum_{m=0}^{\infty}b_{m}\frac{t^m}{m!}
=\sum_{n=k}^{\infty}\binom{n}{k}b_{n-k}\frac{t^{n}}{n!}.
\end{align*}
Thus we get
\begin{align*}
S_{2}(n,k;\bold{P})=\binom{n}{k}b_{n-k},
\end{align*}
where $b_{n}=b_{n}(0)$ are the Bernoulli numbers of the second. \\
Alternatively, it is easily shown that 
\begin{align*}
S_{2}(n,k;\bold{P})=\sum_{m=k}^{n}\sum_{l=k}^{m}\binom{n}{m}b_{n-m}S_{1}(m,l)S_{2}(l,k).
\end{align*}
(s) Let $\bold{P}=\left\{C_{n}(x;a)\right\}$ be the sequence of Poisson-Charlier polynomials. Here $C_{n}(x;a) \sim \big(e^{a(e^{t}-1)}, a(e^{t}-1)\big)$, with $a \ne 0$, and therefore $\sum_{n=0}^{\infty}C_{n}(x;a)\frac{t^{n}}{n!}=e^{-t}\big(1+\frac{t}{a}\big)^{x}$ (see [23]).  From Theorem 3.1, we have
\begin{align*}
\sum_{n=k}^{\infty}S_{2}(n,k;\bold{P})\frac{t^{n}}{n!}&=\frac{1}{k!}e^{-t}\big(e^{\log (1+\frac{t}{a})}-1\big)^{k}=\frac{1}{k!}a^{-k}t^{k}\sum_{m=0}^{\infty}(-1)^{m}\frac{t^{m}}{m!}\\
&=\sum_{n=k}^{\infty}\binom{n}{k}(-1)^{n-k}a^{-k}\frac{t^{n}}{n!}.
\end{align*}
Thus we see that
\begin{align*}
S_{2}(n,k;\bold{P})=\binom{n}{k}(-1)^{n-k}a^{-k}.
\end{align*}
Alternatively, we can also show that
\begin{align*}
S_{2}(n,k;\bold{P})=\sum_{m=k}^{n}\sum_{l=k}^{m}(-1)^{n-m}a^{-m}\binom{n}{m}S_{1}(m,l)S_{2}(l,k).
\end{align*}
(t) Let $\bold{P}=\left\{p_{n}(x)\right\}$, with $p_{n}(x)=\sum_{k=0}^{n}B_{k}(x)B_{n-k}(x)$. This is not a Sheffer sequence. For this, we recall from [13] that
\begin{align*}
p_{n}(x)=\frac{2}{n+2}\sum_{m=0}^{n-2}\binom{n+2}{m}B_{n-m}B_m(x)+(n+1)B_n(x).
\end{align*}
By Theorem 3.1 and (o), we note that
\begin{align*}
\frac{1}{k!}\langle{(e^{t}-1)^{k}|B_{n}(x) \rangle}=\sum_{l=0}^{n-k}\binom{n}{l}S_{2}(n-l,k)B_{l}.
\end{align*}
Again, by Theorem 3.1, we see that
\begin{align*}
&S_{2}(n,k;\bold{P})=\frac{1}{k!}\langle{(e^{t}-1)^{k}|p_{n}(x)\rangle}\\
&=\frac{2}{n+2}\sum_{m=0}^{n-2}\binom{n+2}{m}B_{n-m}\frac{1}{k!}\langle{(e^{t}-1)^{k}|B_{m}(x) \rangle}+(n+1)\frac{1}{k!}\langle{(e^{t}-1)^{k}|B_{n}(x) \rangle}\\
&=\frac{2}{n+2}\sum_{m=0}^{n-2}\binom{n+2}{m}B_{n-m}\sum_{l=0}^{m-k}\binom{m}{l}S_{2}(m-l,k)B_{l}
+(n+1)\sum_{l=0}^{n-k}\binom{n}{l}S_{2}(n-l,k)B_{l}.
\end{align*}
This says that the following identity holds:
\begin{align*}
\sum_{k=0}^{n}B_{k}(x)B_{n-k}(x)&=\sum_{k=0}^{n}\left\{\frac{2}{n+2}\sum_{m=0}^{n-2}\binom{n+2}{m}B_{n-m}\sum_{l=0}^{m-k}\binom{m}{l}S_{2}(m-l,k)B_{l}\right.\\
&\left.\quad\quad\quad+(n+1)\sum_{l=0}^{n-k}\binom{n}{l}S_{2}(n-l,k)B_{l}\right\}(x)_{k}.
\end{align*}

\section{Stirling numbers of the first kind associated with sequences of polynomials}

As in the `classical' case, we would like to introduce the Stirling numbers of the first kind associated with sequences of polynomials. Let $\bold{P}=\left\{p_{n}(x)\right\}_{n=0}^{\infty}$ be a sequence of polynomials such that deg\,$p_{n}(x)=n, p_{0}(x)=1$.
In view of our definition of the Stirling numbers of the second kind associated with $\bold{P}$, it is natural to define the Stirling numbers of the first associated with $\bold{P}$ as the coefficients when we expand $(x)_{n}$ in terms of $p_{k}(x)$:
\begin{align}
(x)_{n}=\sum_{k=0}^{n}S_1(n,k;\bold{P})p_{k}(x).\label{1E}
\end{align}

\begin{theorem}
Let $\bold{P}=\left\{p_{n}(x)\right\}_{n=0}^{\infty}$ be a sequence of polynomials such that deg\,$p_{n}(x)=n, p_{0}(x)=1$. We let
\begin{align*}
(x)_{n}=\sum_{k=0}^{n}S_1(n,k;\bold{P})p_{k}(x).
\end{align*}
(a)\,\, Let $\bar{\bold{P}}=\left\{\bar{p}_{n}(x)\right\}$, where $\bar{p}_{n}(x)$ is the sequence of polynomials defined by $\bar{p}_{0}(x)=1,\,\,\bar{p}_{n}(x)=xp_{n-1}(x),\,\,(n \ge 1)$. Then we have
\begin{align*}
S_{1}(n+1,k;\bar{\bold{P}})=S_{1}(n,k-1;\bold{P})-nS_{1}(n,k;\bar{\bold{P}}), \,\, (0 \le k \le n+1).
\end{align*}
(b)\,\, Assume that $p_{n}(x) \sim (g(t), f(t))$. Then the Stirling numbers of the first kind associated with $\bold{P}$ are given by
\begin{align*}
S_1(n,k;\bold{P})=\frac{1}{k!}\langle{g(t)(f(t))^{k}|(x)_{n}\rangle}=\frac{1}{k!}\sum_{l=k}^{n}S_{1}(n,l)\langle{g(t)f(t)^{k}|x^{l}\rangle}.
\end{align*}
(c)\,\,$(-1)^{n}\sum_{k=0}^{n}S_{1}(n,k;\bold{P})p_{k}(-x)=\langle{x\rangle}_{n}$, \,\,$(-1)^{n}\sum_{k=0}^{n}S_{1}(n,k;\bold{P})p_{k}(-1)=n!.$
\begin{proof}
(a)\,\,This is immediate from the following observation:
\begin{align*}
\sum_{k=0}^{n+1}S_{1}(n+1,k;\bar{\bold{P}})\bar{p}_{k}(x)&=(x)_{n}(x-n)\\
&=x\sum_{k=0}^{n}S_{1}(n,k;\bold{P})p_{k}(x)-n\sum_{k=0}^{n}S_{1}(n,k;\bar{\bold{P}})\bar{p}_{k}(x)\\
&=\sum_{k=0}^{n}S_{1}(n,k;\bold{P})\bar{p}_{k+1}(x)-n\sum_{k=0}^{n}S_{1}(n,k;\bar{\bold{P}})\bar{p}_{k}(x)\\
&=\sum_{k=0}^{n+1}S_{1}(n,k-1;\bold{P})\bar{p}_{k}(x)-n\sum_{k=0}^{n+1}S_{1}(n,k;\bar{\bold{P}})\bar{p}_{k}(x).
\end{align*}
(b)\,\, This follows from \eqref{16B}. \\
(c)\,\, This is clear by definition in \eqref{1E}.
\end{proof}
\end{theorem}
For any delta series $f(t)$, we define the {\it{logarithm associated to $f(t)$}} by
\begin{align*}
L_{f(t)} t =f(\log(1+ t)),
\end{align*}
In case $p_{n}(x) \sim (g(t), f(t))$, $L_{f(t)}t$ may be also denoted by $L_{\bold{P}}t$ and called the {\it{logarithm associated to $\bold{P}=\left\{p_{n}(x)\right\}$}}. Recalling that $(x)_{n} \sim (1,l(t)=e^{t}-1)$, we note that $L_{f(t)}t=f(\bar{l}(t))$.

\begin{theorem}
Let $p_{n}(x) \sim (1, f(t)),$  and let $\bold{P}=\left\{p_{n}(x)\right\}$. Then the generating function of the Stirling numbers of the first kind associated with $\bold{P}$ is given by
\begin{align*}
\sum_{n=k}^{\infty}S_1(n,k;\bold{P})\frac{t^{n}}{n!}=\frac{1}{k!}(L_{f(t)}t)^{k},
\end{align*}
where $L_{f(t)}t=L_{\bold{P}}t=f(\log(1+ t))$ is the logarithm associated to $f(t)$.
\begin{proof}
Assume now that $p_{n}(x) \sim (1, f(t))$. Then we observe that
\begin{align}
\sum_{k=0}^{\infty}p_{k}(x)\sum_{n=k}^{\infty}S_1(n,k;\bold{P})\frac{t^{n}}{n!}&=\sum_{n=0}^{\infty}\sum_{k=0}^{n}S_{1}(n,k;\bold{P})p_{k}(x)\frac{t^{n}}{n!}\label{2E}\\
&=\sum_{n=0}^{\infty}(x)_{n}\frac{t^{n}}{n!}=(1+t)^{x}.\nonumber
\end{align}
On the other hand,  we also have
\begin{align}
\sum_{k=0}^{\infty}p_{k}(x)\frac{1}{k!}\big(L_{f(t)}t \big)^{k}&=e^{x\bar{f}\big(f (\log(1+t))\big)}=(1+t)^{x},\label{3E}
\end{align}
where we used the fact that $\sum_{k=0}^{\infty}p_{k}(x)\frac{t^{k}}{k!}=e^{x \bar{f}(t)}$.  Now,  the assertion follows from \eqref{2E} and \eqref{3E}. \\
We would like to give another proof for this by using Theorem 5.1 (b). \\
Here we have to show that
\begin{align*}
\sum_{n=k}^{\infty}\frac{1}{k!}\langle{f(t)^{k}|(x)_{n}\rangle}\frac{t^{n}}{n!}=\frac{1}{k!}(f(\log(1+t)))^{k}.
\end{align*}
Equivalently, we need to verify that
\begin{align}
\sum_{n=k}^{\infty}\langle{f(t)^{k}|(x)_{n}\rangle}\frac{(e^{t}-1)^{n}}{n!}=(f(t))^{k}. \label{4E}
\end{align}
The left hand side of \eqref{4E} is equal to
\begin{align*}
\sum_{n=k}^{\infty}\sum_{l=k}^{n}S_{1}(n,l)\langle{f(t)^{k}|x^{l}\rangle}\frac{(e^{t}-1)^{n}}{n!}&=\sum_{l=k}^{\infty}\langle{f(t)^{k}|x^{l}\rangle}\sum_{n=l}^{\infty}S_{1}(n,l)\frac{(e^{t}-1)^{n}}{n!}\\
&=\sum_{l=0}^{\infty}\langle{f(t)^{k}|x^{l}\rangle}\frac{t^{l}}{l!}=(f(t))^{k}.
\end{align*}
\end{proof}
\end{theorem}

The following orthogonality and inverse relations can be easily shown.
\begin{theorem}
Let $\bold{P}=\left\{p_{n}(x)\right\}_{n=0}^{\infty}$ be a sequence of polynomials such that deg\,$p_{n}(x)=n, p_{0}(x)=1$. Then we have
\begin{align*}
&(a)\,\, \sum_{k=l}^{n}S_{1}(n,k;\bold{P})S_{2}(k,l;\bold{P})=\delta_{n,l},
\,\, \sum_{k=l}^{n}S_{2}(n,k;\bold{P})S_{1}(k,l;\bold{P})=\delta_{n,l},\,\,(0 \le l \le n). \\
&(b)\,\, a_{n}=\sum_{k=0}^{n}S_{2}(n,k;\bold{P})c_{k} \Longleftrightarrow c_{n}=\sum_{k=0}^{n}S_{1}(n,k;\bold{P})a_{k}.\\
&(c)\,\, a_{n}=\sum_{k=n}^{m}S_{2}(k,n;\bold{P})c_{k} \Longleftrightarrow c_{n}=\sum_{k=n}^{m}S_{1}(k,n;\bold{P})a_{k}.
\end{align*}
\end{theorem}

\begin{remark}
Let $p_{n}(x)$ be Sheffer for the pair $(g(t), f(t))$. Then we saw in Theorem 3.2 that the generating function of $S_{2}(n,k;\bold{P})$ is given by
\begin{align*}
\sum_{n=k}^{\infty}S_2(n,k;\bold{P})\frac{t^{n}}{n!}=\frac{1}{g(\bar{f}(t))}\frac{1}{k!}(e^{\bar{f}(t)}-1)^{k}.
\end{align*}
Let $f(t)$ be any delta series. After introducing the logarithm associated to $f(t)$, it seems natural to introduce the {\it{exponential associated to $f(t)$}} as 
\begin{align*}
E_{f(t)}t=e^{\bar{f}(t)}-1.
\end{align*}
When $p_{n}(x) \sim (g(t), f(t))$, $E_{f(t)}t$ may be also denoted by $E_{\bold{P}}t$ and called the exponential associated to $\bold{P}$. Recalling that $(x)_{n} \sim (1,l(t)=e^{t}-1)$, we observe that $E_{f(t)}t=l(\bar{f}(t))$. In summary, we have
\begin{align*}
\sum_{n=k}^{\infty}S_2(n,k;\bold{P})\frac{t^{n}}{n!}=\frac{1}{g(\bar{f}(t))}\frac{1}{k!}(E_{f(t)}t)^{k}.\,
\end{align*}
where $E_{f(t)}t=E_{\bold{P}}t=e^{\bar{f}(t)}-1$.\\
As an example, consider the delta series $ f(t)=\frac{1}{\lambda}(e^{\lambda t}-1)$. Here $\bar{f}(t)=\frac{1}{\lambda}\log(1+\lambda t)$. Thus we see that
\begin{align*}
L_{f(t)}t=\frac{1}{\lambda}\big((1+t)^{\lambda}-1 \big),\quad E_{f(t)}t=(1+\lambda t)^{\frac{1}{\lambda}}-1.
\end{align*}
Here we note that, in many of the papers on degenerate special polynomials and numbers, $\frac{1}{\lambda}\big(t^{\lambda}-1 \big)$ is called the degenerate logarithm and denoted by $\log_{\lambda}t$, and $(1+\lambda t)^{\frac{1}{\lambda}}$ is called the degenerate exponential and denoted by $e_{\lambda}(t)$.
\end{remark}

\section{Examples on Stirling numbers of the first kind associated with sequences of polynomials}
(a) Let $\bold{P}=\left\{x^{n}\right\}$. Then $x^{n} \sim (1,t)$, and $S_{2}(n,k;\bold{P})=S_{2}(n,k)$.  As $(x)_{n}=\sum_{k=0}^{n}S_{1}(n,k;\bold{P})x^{k}$, we have (see [4,22,23])
\begin{align*}
S_{1}(n,k;\bold{P})=S_{1}(n,k), \quad \sum_{n=k}^{\infty}S_{1}(n,k;\bold{P})\frac{t^n}{n!}=\frac{1}{k!}(\log(1+t))^{k}.
\end{align*}
(b) Let $\bold{P}=\left\{(x)_{n,\lambda}\right\}$ be the sequence of generalized falling factorials. Then $(x)_{n,\lambda} \sim \big(1, f(t)=\frac{1}{\lambda}(e^{\lambda t}-1)\big)$, and $S_{2}(n,k;\bold{P})=S_{2,\lambda}(n,k)$. As $(x)_{n}=\sum_{k=0}^{n}S_{1}(n,k;\bold{P})(x)_{k,\lambda}$, we have (see [9,10,12,14,16-18])
\begin{align*}
S_{1}(n,k;\bold{P})=S_{1,\lambda}(n,k),\quad \sum_{n=k}^{\infty}S_{1,\lambda}(n,k)\frac{t^n}{n!}=\frac{1}{k!}(\log_{\lambda}(1+t))^{k},
\end{align*}
(see Remark 5.4).\\ Now, by Theorem 5.3, we get
\begin{align*}
&\sum_{k=l}^{n}S_{1,\lambda}(n,k)S_{2,\lambda}(k,l)=\delta_{n,l},
\,\, \sum_{k=l}^{n}S_{2,\lambda}(n,k)S_{1,\lambda}(k,l)=\delta_{n,l},\,\,(0 \le l \le n), \\
&a_{n}=\sum_{k=0}^{n}S_{2,\lambda}(n,k)c_{k} \Longleftrightarrow c_{n}=\sum_{k=0}^{n}S_{1,\lambda}(n,k)a_{k},\\
&a_{n}=\sum_{k=n}^{m}S_{2,\lambda}(k,n)c_{k} \Longleftrightarrow c_{n}=\sum_{k=n}^{m}S_{1,\lambda}(k,n)a_{k}.
\end{align*}
(c) Let $\bold{P}=\left\{\langle{x \rangle}_{n}\right\}$ be the sequence of rising factorials. Then $\langle{x \rangle}_{n} \sim (1, 1-e^{-t})$, and $S_{2}(n,k;\bold{P})=L(n,k)$ are the (unsigned) Lah numbers. In addition, by Theorem 5.2, we have \\
\begin{align*}
\sum_{n=k}^{\infty}S_1(n,k;\bold{P})\frac{t^{n}}{n!}=\frac{1}{k!}\bigg(1-\frac{1}{1+t}\bigg)^{k}=\sum_{n=k}^{\infty}(-1)^{n-k}L(n,k)\frac{t^{n}}{n!}.
\end{align*}
Thus $S_{1}(n,k;\bold{P})=(-1)^{n-k}L(n,k)$ are the signed Lah numbers (see [4,11,22,23]).
Hence,  by Theorem 5.3,  we have the identities about Lah numbers given by:
\begin{align}
&\sum_{k=l}^{n}(-1)^{n-k}L(n,k)L(k,l)=\delta_{n,l},\quad \sum_{k=l}^{n}(-1)^{k-l}L(n,k)L(k,l)=\delta_{n,l}, \label{1F}\\
&a_{n}=\sum_{k=0}^{n}L(n,k)c_{k} \Longleftrightarrow c_{n}=\sum_{k=0}^{n}(-1)^{n-k}L(n,k)a_{k},\nonumber \\
&a_{n}=\sum_{k=n}^{m}L(k,n)c_{k} \Longleftrightarrow c_{n}=\sum_{k=n}^{m}(-1)^{k-n}L(k,n)a_{k}. \nonumber
\end{align}
Note here that the two identities in \eqref{1F} are equivalent. \\
(d) Let $\bold{P}=\left\{\langle{x \rangle}_{n,\lambda}\right\}$ be the sequence of generalized rising factorials. Then $\langle{x \rangle}_{n,\lambda} \sim \big(1, \frac{1}{\lambda}(1-e^{-\lambda t})\big)$, and $S_{2}(n,k;\bold{P})=L_{\lambda}(n,k)$ are the degenerate Lah numbers given by 
\begin{align*}
\sum_{n=k}^{\infty}L_{\lambda}(n,k)\frac{t^{n}}{n!}&=\frac{1}{k!}\big(e^{-\frac{1}{\lambda}\log (1-\lambda t)}-1\big)^{k}=\frac{1}{k!}\big((1-\lambda t)^{-\frac{1}{\lambda}}-1\big)^{k}.
\end{align*}
By Theorem 5.2, we have
\begin{align*}
\sum_{n=k}^{\infty}S_{1}(n,k;\bold{P})\frac{t^{n}}{n!}&=\frac{1}{k!}\big(\frac{1}{\lambda}(1-(1+t)^{-\lambda})\big)^{k} \\
&=\sum_{n=k}^{\infty}(-\lambda)^{n-k}L_{\frac{1}{\lambda}}(n,k)\frac{t^{n}}{n!}.
\end{align*}
 Thus $S_{1}(n,k;\bold{P})=(-\lambda)^{n-k}L_{\frac{1}{\lambda}}(n,k)$ (see [14]).  Hence,  by Theorem 5.3,  we get
\begin{align*}
&\sum_{k=l}^{n}(-\lambda)^{n-k}L_{\frac{1}{\lambda}}(n,k)L_{\lambda}(k,l)=\delta_{n,l},\,\,
\sum_{k=l}^{n}(-\lambda)^{k-l}L_{\lambda}(n,k)L_{\frac{1}{\lambda}}(k,l)=\delta_{n,l},\\
&a_{n}=\sum_{k=0}^{n}L_{\lambda}(n,k)c_{k} \Longleftrightarrow c_{n}=\sum_{k=0}^{n}(-\lambda)^{n-k}L_{\frac{1}{\lambda}}(n,k)a_{k},\\
&a_{n}=\sum_{k=n}^{m}L_{\lambda}(k,n)c_{k} \Longleftrightarrow c_{n}=\sum_{k=n}^{m}(-\lambda)^{k-n}L_{\frac{1}{\lambda}}(k,n)a_{k}.
\end{align*}
We remark that $S_{1}(n,k;\bold{P})=(-\lambda)^{n-k}L_{\frac{1}{\lambda}}(n,k)=S_{1,-\lambda}(n,k).$ Indeed,  we see that\\
\begin{align*}
\sum_{n=k}^{\infty}S_{1}(n,k;\bold{P})\frac{t^{n}}{n!}&=\frac{1}{k!}\big(\frac{1}{-\lambda}((1+t)^{-\lambda})-1\big)^{k} \\
&=\frac{1}{k!}\big(\log_{-\lambda}(1+t)\big)^{k}\\
&=\sum_{n=k}^{\infty}S_{1,-\lambda}(n,k)\frac{t^{n}}{n!}.
\end{align*}
(e) Let $\bold{P}=\left\{x^{[n]}\right\}$ be the sequence of central factorials. Then $x^{[n]} \sim (1, f(t)=e^{\frac{t}{2}}-e^{-\frac{t}{2}})$, and $S_{2}(n,k;\bold{P})=\sum_{l=k}^{n}T_{1}(n,l)S_{2}(l,k)$, 
where $T_{1}(n,k)$ are the central factorial numbers of the first kind defined by $x^{[n]}=\sum_{k=0}^{n}T_{1}(n,k)x^{k}$. Recall here that
\begin{align*}
\frac{1}{k!}\big(e^{\frac{t}{2}}-e^{-\frac{t}{2}}\big)^{k}=\sum_{n=k}^{\infty}T_{2}(n,k)\frac{t^{n}}{n!}.
\end{align*} \\
Now,  by Theorem 5.2, we obtain
\begin{align*}
\sum_{n=k}^{\infty}S_{1}(n,k;\bold{P})\frac{t^{n}}{n!}&=\sum_{l=k}^{\infty}T_{2}(l,k)\frac{1}{l!}(\log(1+t))^{l}\\
&=\sum_{n=k}^{\infty}\sum_{l=k}^{n}S_{1}(n,l)T_{2}(l,k)\frac{t^{n}}{n!}.
\end{align*}
Hence $S_{1}(n,k;\bold{P})=\sum_{l=k}^{n}S_{1}(n,l)T_{2}(l,k)$ (see [3,4,15,23]).   From Theorem 5.3, we get:
\begin{align*}
&\sum_{k=l}^{n}\sum_{m=k}^{n}\sum_{j=l}^{k}S_{1}(n,m)T_{2}(m,k)T_{1}(k,j)S_{2}(j,l)=\delta_{n,l}, \\
&\sum_{k=l}^{n}\sum_{j=k}^{n}\sum_{m=l}^{k}T_{1}(n,j)S_{2}(j,k)S_{1}(k,m)T_{2}(m,l)=\delta_{n,l},\\
&a_{n}=\sum_{k=0}^{n}\sum_{l=k}^{n}T_{1}(n,l)S_{2}(l,k)c_{k} \Longleftrightarrow c_{n}=\sum_{k=0}^{n}\sum_{l=k}^{n}S_{1}(n,l)T_{2}(l,k)a_{k},\\
&a_{n}=\sum_{k=n}^{m}\sum_{l=n}^{k}T_{1}(k,l)S_{2}(l,n)c_{k} \Longleftrightarrow c_{n}=\sum_{k=n}^{m}\sum_{l=n}^{k}S_{1}(k,l)T_{2}(l,n)a_{k}.
\end{align*}
(f) Let $\bold{P}=\left\{\mathrm{Bel}_{n}^{(c)}(x)\right\}$ be the sequence of central Bell polynomials. Then $\mathrm{Bel}_{n}^{(c)}(x) \sim \Big(1, f(t)=2\log\big(\frac{t+\sqrt{t^{2}+4}}{2}\big)\Big)$, and $S_{2}(n,k;\bold{P})=\sum_{l=k}^{n}T_{2}(n,l)S_{2}(l,k),$ where $T_{2}(n,k)$ are the central factorial numbers of the second kind defined by $x^{n}=\sum_{k=0}^{n}T_{2}(n,k)x^{[k]}$.
We recall here that
\begin{align*}
\frac{1}{k!}\bigg(2\log\bigg(\frac{t+\sqrt{t^{2}+4}}{2}\bigg)\bigg)^{k}
=\sum_{n=k}^{\infty}T_{1}(n,k)\frac{t^{n}}{n!}.
\end{align*}
Now, by Theorem 5.2, we get
\begin{align*}
\sum_{n=k}^{\infty}S_{1}(n,k;\bold{P})\frac{t^{n}}{n!}&=\sum_{l=k}^{\infty}T_{1}(l,k)\frac{1}{l!}(\log(1+t))^{l}\\
&=\sum_{n=k}^{\infty}\sum_{l=k}^{n}S_{1}(n,l)T_{1}(l,k)\frac{t^{n}}{n!}.
\end{align*}
Therefore $S_{1}(n,k;\bold{P})=\sum_{l=k}^{n}S_{1}(n,l)T_{1}(l,k)$ (see [15]).  From Theorem 5.3,  we have:
\begin{align*}
&\sum_{k=l}^{n}\sum_{m=k}^{n}\sum_{j=l}^{k}S_{1}(n,m)T_{1}(m,k)T_{2}(k,j)S_{2}(j,l)=\delta_{n,l},\\
&\sum_{k=l}^{n}\sum_{j=k}^{n}\sum_{m=l}^{k}T_{2}(n,j)S_{2}(j,k)S_{1}(k,m)T_{1}(m,l)=\delta_{n,l},\\
&a_{n}=\sum_{k=0}^{n}\sum_{l=k}^{n}T_{2}(n,l)S_{2}(l,k)c_{k} \Longleftrightarrow c_{n}=\sum_{k=0}^{n}\sum_{l=k}^{n}S_{1}(n,l)T_{1}(l,k)a_{k},\\
&a_{n}=\sum_{k=n}^{m}\sum_{l=n}^{k}T_{2}(k,l)S_{2}(l,n)c_{k} \Longleftrightarrow c_{n}=\sum_{k=n}^{m}\sum_{l=n}^{k}S_{1}(k,l)T_{1}(l,n)a_{k}.
\end{align*}

(g) Let $\bold{P}=\left\{\mathrm{Bel}_{n,\lambda}^{(c)}(x)\right\}$ be the sequence of degenerate central Bell polynomials. Then $\mathrm{Bel}_{n,\lambda}^{(c)}(x) \sim \Big(1, f(t)=\frac{1}{\lambda}\Big(\Big(\frac{t+\sqrt{t^{2}+4}}{2}\Big)^{2 \lambda}-1\Big)\Big),$ and $S_{2}(n,k;\bold{P})=\sum_{l=k}^{n}T_{2,\lambda}(n,l)S_{2}(l,k)$, where $T_{2,\lambda}(n,k)$ are the degenerate central factorial numbers of the second kind defined by $(x)_{n,\lambda}=\sum_{k=0}^{n}T_{2,\lambda}(n,k)x^{[k]}$. Here we recall that
\begin{align*}
\frac{1}{k!}\left\{\frac{1}{\lambda}\bigg(\bigg(\frac{t+\sqrt{t^{2}+4}}{2}\bigg)^{2 \lambda}-1\bigg)\right\}^{k}=\sum_{n=k}^{\infty}T_{1,\lambda}(n,k)\frac{t^{n}}{n!}.
\end{align*}
Now, by Theorem 5.2, we have
\begin{align*}
\sum_{n=k}^{\infty}S_{1}(n,k;\bold{P})\frac{t^{n}}{n!}&=\sum_{l=k}^{\infty}T_{1,\lambda}(l,k)\frac{1}{l!}(\log(1+t))^{l}\\
&=\sum_{n=k}^{\infty}\sum_{l=k}^{n}S_{1}(n,l)T_{1,\lambda}(l,k)\frac{t^{n}}{n!}.
\end{align*}
Thus we have $S_{1}(n,k;\bold{P})=\sum_{l=k}^{n}S_{1}(n,l)T_{1,\lambda}(l,k)$ (see [16]).
From Theorem 5.3, we have:
\begin{align*}
&\sum_{k=l}^{n}\sum_{m=k}^{n}\sum_{j=l}^{k}S_{1}(n,m)T_{1,\lambda}(m,k)T_{2,\lambda}(k,j)S_{2}(j,l)=\delta_{n,l},\\
&\sum_{k=l}^{n}\sum_{j=k}^{n}\sum_{m=l}^{k}T_{2,\lambda}(n,j)S_{2}(j,k)S_{1}(k,m)T_{1,\lambda}(m,l)=\delta_{n,l},\\
&a_{n}=\sum_{k=0}^{n}\sum_{l=k}^{n}T_{2,\lambda}(n,l)S_{2}(l,k)c_{k} \Longleftrightarrow c_{n}=\sum_{k=0}^{n}\sum_{l=k}^{n}S_{1}(n,l)T_{1,\lambda}(l,k)a_{k},\\
&a_{n}=\sum_{k=n}^{m}\sum_{l=n}^{k}T_{2,\lambda}(k,l)S_{2}(l,n)c_{k} \Longleftrightarrow c_{n}=\sum_{k=n}^{m}\sum_{l=n}^{k}S_{1}(k,l)T_{1,\lambda}(l,n)a_{k}.
\end{align*}

(h) Let $x^{[n,\lambda]}$ be the sequence of polynomials defined by $x^{[0,\lambda]}=1,\quad x^{[n,\lambda]}=x\big(x+(\frac{1}{2}n-1)\lambda\big)_{n-1,\lambda}, \,\,(n \ge 1)$. Then $x^{[n,\lambda]} \sim \big(1, f(t)=\frac{1}{\lambda}(e^{\frac{\lambda t}{2}}-e^{-\frac{\lambda t}{2}})\big)$, and $S_{2}(n,k;\bold{P})=\sum_{l=k}^{n}R_{1,\lambda}(n,l)S_{2}(l,k)$,
where $R_{1,\lambda}(n,k)$ are defined by $x^{[n,\lambda]}=\sum_{k=0}^{n}R_{1,\lambda}(n,k)x^{k}$. Here we recall that
\begin{align*}
 \sum_{n=k}^{\infty}R_{2,\lambda}(n,k)\frac{t^{n}}{n!}=\frac{1}{k!}\Big(\frac{1}{\lambda}\big(e^{\frac{\lambda t}{2}}-e^{-\frac{\lambda t}{2}}\big)\Big)^{k}.
\end{align*}
Let $\bold{P}=\left\{x^{[n,\lambda]}\right\}$. Now, by Theorem 5.2, we have
\begin{align*}
\sum_{n=k}^{\infty}S_{1}(n,k;\bold{P})\frac{t^{n}}{n!}&=\sum_{l=k}^{\infty}R_{2,\lambda}(l,k)\frac{1}{l!}\big(\log (1+t)\big)^{l}\\
&=\sum_{n=k}^{\infty}\sum_{l=k}^{n}S_{1}(n,l)R_{2,\lambda}(l,k)\frac{t^{n}}{n!}.
\end{align*}
Hence $S_{1}(n,k;\bold{P})=\sum_{l=k}^{n}S_{1}(n,l)R_{2,\lambda}(l,k)$. From Theorem 5.3, we get
\begin{align*}
&\sum_{k=l}^{n}\sum_{m=k}^{n}\sum_{j=l}^{k}S_{1}(n,m)R_{2,\lambda}(m,k)R_{1,\lambda}(k,j)S_{2}(j,l)=\delta_{n,l}, \\
&\sum_{k=l}^{n}\sum_{j=k}^{n}\sum_{m=l}^{k}R_{1,\lambda}(n,j)S_{2}(j,k)S_{1}(k,m)R_{2,\lambda}(m,l)=\delta_{n,l},\\
&a_{n}=\sum_{k=0}^{n}\sum_{l=k}^{n}R_{1,\lambda}(n,l)S_{2}(l,k)c_{k} \Longleftrightarrow c_{n}=\sum_{k=0}^{n}\sum_{l=k}^{n}S_{1}(n,l)R_{2,\lambda}(l,k)a_{k},\\
&a_{n}=\sum_{k=n}^{m}\sum_{l=n}^{k}R_{1,\lambda}(k,l)S_{2}(l,n)c_{k} \Longleftrightarrow c_{n}=\sum_{k=n}^{m}\sum_{l=n}^{k}S_{1}(k,l)R_{2,\lambda}(l,n)a_{k}.
\end{align*}

(i) Let $\bold{P}=\left\{B_{n}^{L}(x)\right\}$ be the sequence of Lah-Bell polynomials. Then $B_{n}^{L}(x) \sim (1,\frac{t}{1+t})$, and $S_{2}(n,k;\bold{P})=\sum_{l=k}^{n}L(n,l)S_{2}(l,k)$. \\
By Theorem 5.1 and (c) above,  we see that
\begin{align*}
S_{1}(n,k;\bold{P})&=\sum_{l=k}^{n} S_{1}(n,l)\big\langle{\frac{1}{k!}\Big(\frac{t}{1+t}\Big)^{k}|x^{l} \big\rangle}\\
&=\sum_{l=k}^{n}S_{1}(n,l)\sum_{j=k}^{l}(-1)^{j-k}L(j,k)\frac{1}{j!}\langle{t^{j}|x^{l}\rangle}\\
&=\sum_{l=k}^{n}(-1)^{l-k}S_{1}(n,l)L(l,k).
\end{align*}
Therefore we have $S_{1}(n,k;\bold{P})=\sum_{l=k}^{n}(-1)^{l-k}S_{1}(n,l)L(l,k)$ (see [11]).  Now,  from Theorem 5.3, we obtain 
\begin{align*}
&\sum_{k=l}^{n}\sum_{m=k}^{n}\sum_{j=l}^{k}(-1)^{m-k}S_{1}(n,m)L(m,k)L(k,j)S_2(j,l)=\delta_{n,l}, \\
&\sum_{k=l}^{n}\sum_{j=k}^{n}\sum_{m=l}^{k}(-1)^{m-l}L(n,j)S_{2}(j,k)S_{1}(k,m)L(m,l)=\delta_{n,l},\\
&a_{n}=\sum_{k=0}^{n}\sum_{l=k}^{n}L(n,l)S_{2}(l,k)c_{k} \Longleftrightarrow c_{n}=\sum_{k=0}^{n}\sum_{l=k}^{n}(-1)^{l-k}S_{1}(n,l)L(l,k)a_{k},\\
&a_{n}=\sum_{k=n}^{m}\sum_{l=n}^{k}L(k,l)S_{2}(l,n)c_{k} \Longleftrightarrow c_{n}=\sum_{k=n}^{l}\sum_{l=n}^{k}(-1)^{l-n}S_{1}(k,l)L(l,n)a_{k}.
\end{align*}
(j) Let $\bold{P}=\left\{B_{n,\lambda}^{L}(x)\right\}$ be the sequence of degenerate Lah-Bell polynomials. Then $B_{n,\lambda}^{L}(x) \sim (1,\frac{e^{\lambda t}-1}{\lambda+e^{\lambda t}-1})$, and $S_{2}(n,k;\bold{P})=\sum_{l=k}^{n}L(n,l)L_{-\lambda}(l,k).$ Now, by Theorem 5.2, we get
\begin{align*}
\sum_{n=k}^{\infty}S_{1}(n,k;\bold{P})\frac{t^{n}}{n!}&=\frac{1}{k!}\bigg(\frac{\log_{\lambda}(1+t)}{1+\log_{\lambda}(1+t)}\bigg)^{k}\\
&=\sum_{l=k}^{\infty}(-1)^{l-k}L(l,k)\frac{1}{l!}(\log_{\lambda}(1+t))^{l}\\
&=\sum_{n=k}^{\infty}\sum_{l=k}^{n}(-1)^{l-k}S_{1,\lambda}(n,l)L(l,k)\frac{t^{n}}{n!}.
\end{align*}
Thus we have $S_{1}(n,k;\bold{P})=\sum_{l=k}^{n}(-1)^{l-k}S_{1,\lambda}(n,l)L(l,k)$ (see [14]).
Now, by Theorem 5.3, we obtain
\begin{align*}
&\sum_{k=l}^{n}\sum_{m=k}^{n}\sum_{j=l}^{k}(-1)^{m-k}S_{1,\lambda}(n,m)L(m,k)L(k,j)L_{-\lambda}(j,l)=\delta_{n,l},\\
&\sum_{k=l}^{n}\sum_{j=k}^{n}\sum_{m=l}^{k}(-1)^{m-l}L(n,j)L_{-\lambda}(j,k)S_{1,\lambda}(k,m)L(m,l)=\delta_{n,l},\\
&a_{n}=\sum_{k=0}^{n}\sum_{l=k}^{n}L(n,l)L_{-\lambda}(l,k)c_{k} \Longleftrightarrow c_{n}=\sum_{k=0}^{n}\sum_{l=k}^{n}(-1)^{l-k}S_{1,\lambda}(n,l)L(l,k)a_{k},\\
&a_{n}=\sum_{k=n}^{m}\sum_{l=n}^{k}L(k,l)L_{-\lambda}(l,n)c_{k} \Longleftrightarrow c_{n}=\sum_{k=n}^{m}\sum_{l=n}^{k}(-1)^{l-n}S_{1,\lambda}(k,l)L(l,n)a_{k}.
\end{align*}
(k) Let $\bold{P}=\left\{\mathrm{Bel}_{n}(x)\right\}$ be the sequence of Bell polynomials. Then $\mathrm{Bel}_{n}(x) \sim (1, \log (1+t))$, and $S_{2}(n,k;\bold{P})=\sum_{l=k}^{n}S_2(n,l)S_{2}(l,k)$.
By Theorem 5.1, we get
\begin{align*}
S_{1}(n,k;\bold{P})&=\sum_{l=k}^{n}S_{1}(n,l)\langle{\frac{1}{k!}(\log(1+t))^{k}|x^{l}\rangle}\\
&=\sum_{l=k}^{n}S_{1}(n,l)\sum_{j=k}^{l}S_{1}(j,k)\frac{1}{j!}\langle{t^{j}|x^{l}\rangle}\\
&=\sum_{l=k}^{n}S_{1}(n,l)S_{1}(l,k).
\end{align*}
Thus we see that $S_{1}(n,k;\bold{P})=\sum_{l=k}^{n}S_{1}(n,l)S_{1}(l,k)$ (see [4,22,23]).  Now,  by Theorem 5.3,  we have
\begin{align*}
&\sum_{k=l}^{n}\sum_{m=k}^{n}\sum_{j=l}^{k}S_{1}(n,m)S_{1}(m,k)S_{2}(k,j)S_{2}(j,l)=\delta_{n,l}, \\
&\sum_{k=l}^{n}\sum_{j=k}^{n}\sum_{m=l}^{k}S_{2}(n,j)S_{2}(j,k)S_{1}(k,m)S_{1}(m,l)=\delta_{n,l},\\
&a_{n}=\sum_{k=0}^{n}\sum_{l=k}^{n}S_2(n,l)S_{2}(l,k)c_{k} \Longleftrightarrow c_{n}=\sum_{k=0}^{n}\sum_{l=k}^{n}S_{1}(n,l)S_{1}(l,k)a_{k},\\
&a_{n}=\sum_{k=n}^{m}\sum_{l=n}^{k}S_2(k,l)S_{2}(l,n)c_{k} \Longleftrightarrow c_{n}=\sum_{k=n}^{m}\sum_{l=n}^{k}S_{1}(k,l)S_{1}(l,n)a_{k}.
\end{align*}

(l) Let $\bold{P}=\left\{\mathrm{Bel}_{n,\lambda}(x)\right\}$ be the sequence of partially degenerate Bell polynomials. Then $\mathrm{Bel}_{n,\lambda}(x) \sim (1, \log_{\lambda}(1+t))$, and 
$S_{2}(n,k;\bold{P})=\sum_{l=k}^{n}S_{2,\lambda}(n,l)S_{2}(l,k)$. By Theorem 5.1,  we have
\begin{align*}
S_{1}(n,k;\bold{P})&=\sum_{l=k}^{n}S_{1}(n,l)\langle{\frac{1}{k!}(\log_{\lambda}(1+t))^{k}|x^{l}\rangle}\\
&=\sum_{l=k}^{n}S_{1}(n,l)\sum_{j=k}^{l}S_{1,\lambda}(j,k)\frac{1}{j!}\langle{t^{j}|x^{l}\rangle}\\
&=\sum_{l=k}^{n}S_{1}(n,l)S_{1,\lambda}(l,k).
\end{align*}
Hence we have $S_{1}(n,k;\bold{P})=\sum_{l=k}^{n}S_{1}(n,l)S_{1,\lambda}(l,k)$ (see [17]).  Now,  by Theorem 5.3,  we obtain 
\begin{align*}
&\sum_{k=l}^{n}\sum_{m=k}^{n}\sum_{j=l}^{k}S_{1}(n,m)S_{1,\lambda}(m,k)S_{2,\lambda}(k,j)S_{2}(j,l)=\delta_{n,l}, \\
&\sum_{k=l}^{n}\sum_{j=k}^{n}\sum_{m=l}^{k}S_{2,\lambda}(n,j)S_{2}(j,k)S_{1}(k,m)S_{1,\lambda}(m,l)=\delta_{n,l},
\end{align*}
\begin{align*}
&a_{n}=\sum_{k=0}^{n}\sum_{l=k}^{n}S_{2,\lambda}(n,l)S_{2}(l,k)c_{k} \Longleftrightarrow c_{n}=\sum_{k=0}^{n}\sum_{l=k}^{n}S_{1}(n,l)S_{1,\lambda}(l,k)a_{k},\\
&a_{n}=\sum_{k=n}^{m}\sum_{l=n}^{k}S_{2,\lambda}(k,l)S_{2}(l,n)c_{k} \Longleftrightarrow c_{n}=\sum_{k=n}^{m}\sum_{l=n}^{k}S_{1}(k,l)S_{1,\lambda}(l,n)a_{k}.
\end{align*}
Let $\bold{P}=\left\{\phi_{n,\lambda}(x)\right\}$ be the sequence of fully degenerate Bell polynomials. Then $\phi_{n,\lambda}(x) \sim \big(1, \log_{\lambda}(1+\frac{1}{\lambda}(e^{\lambda t}-1))\big)$, and $S_{2}(n,k;\bold{P})=\sum_{l=k}^{n}S_{2,\lambda}(n,l)L_{-\lambda}(l,k).$ By Theorem 5.2, we obtain
\begin{align*}
\sum_{n=k}^{\infty}S_{1}(n,k;\bold{P})\frac{t^{n}}{n!}&=\frac{1}{k!}\big(\log_{\lambda}(1+\log_{\lambda}(1+t))\big)^{k}=\sum_{l=k}^{\infty}S_{1,\lambda}(l,k)\frac{1}{l!}\big(\log_{\lambda}(1+t)\big)^{l}\\
&=\sum_{l=k}^{\infty}S_{1,\lambda}(l,k)\sum_{n=l}^{\infty}S_{1,\lambda}(n,l)\frac{t^{n}}{n!}
=\sum_{n=k}^{\infty}\sum_{l=k}^{n}S_{1,\lambda}(n,l)S_{1,\lambda}(l,k)\frac{t^{n}}{n!}.
\end{align*}
Therefore we get $S_{1}(n,k;\bold{P})=\sum_{l=k}^{n}S_{1,\lambda}(n,l)S_{1,\lambda}(l,k)$ (see [18]).
Now, by Theorem 5.3, we obtain 
\begin{align*}
&\sum_{k=l}^{n}\sum_{m=k}^{n}\sum_{j=l}^{k}S_{1,\lambda}(n,m)S_{1,\lambda}(m,k)S_{2,\lambda}(k,j)L_{-\lambda}(j,l)=\delta_{n,l}, \\
&\sum_{k=l}^{n}\sum_{j=k}^{n}\sum_{m=l}^{k}S_{2,\lambda}(n,j)L_{-\lambda}(j,k)S_{1,\lambda}(k,m)S_{1,\lambda}(m,l)=\delta_{n,l},\\
&a_{n}=\sum_{k=0}^{n}\sum_{l=k}^{n}S_{2,\lambda}(n,l)L_{-\lambda}(l,k)c_{k} \Longleftrightarrow c_{n}=\sum_{k=0}^{n}\sum_{l=k}^{n}S_{1,\lambda}(n,l)S_{1,\lambda}(l,k)a_{k},\\
&a_{n}=\sum_{k=n}^{m}\sum_{l=n}^{k}S_{2,\lambda}(k,l)L_{-\lambda}(l,n)c_{k} \Longleftrightarrow c_{n}=\sum_{k=n}^{m}\sum_{l=n}^{k}S_{1,\lambda}(k,l)S_{1,\lambda}(l,n)a_{k}.
\end{align*}
(m) Let $\bold{P}=\left\{M_{n}(x)\right\}$ be the sequence of Mittag-Leffler polynomials. Then $M_{n}(x) \sim (1,f(t)=\frac{e^{t}-1}{e^{t}+1})$, and $S_{2}(n,k;\bold{P})=2^{k}L(n,k)$. Then, from Theorem 5.2, we have
\begin{align*}
\sum_{n=k}^{\infty}S_{1}(n,k;\bold{P})\frac{t^{n}}{n!}&=\frac{1}{k!}\Big(\frac{\frac{t}{2}}{1+\frac{t}{2}}\Big)^{k}=\sum_{n=k}^{\infty}(-1)^{n-k}L(n,k)\frac{1}{2^{n}}\frac{t^{n}}{n!}.
\end{align*}
Hence we have $S_{1}(n,k;\bold{P})=(-1)^{n-k}L(n,k)\frac{1}{2^{n}}$ (see [23]).  By Theorem 5.3, we get
\begin{align*}
&\sum_{k=l}^{n}(-1)^{n-k}L(n,k)L(k,l)=\delta_{n,l}, \\
&\sum_{k=l}^{n}(-1)^{k-l}L(n,k)L(k,l)=\delta_{n,l},\\
&a_{n}=\sum_{k=0}^{n}2^{k}L(n,k)c_{k} \Longleftrightarrow c_{n}=\sum_{k=0}^{n}(-1)^{n-k}L(n,k)\frac{1}{2^{n}}a_{k},\\
&a_{n}=\sum_{k=n}^{m}2^{n}L(k,n)c_{k} \Longleftrightarrow c_{n}=\sum_{k=n}^{m}(-1)^{k-n}L(k,n)\frac{1}{2^{k}}a_{k}.
\end{align*}
Note that the same identities were obtained in (c).

(n) Let $\bold{P}=\left\{L_{n}(x)\right\}$ be the sequence of Laguerre polynomials of order -1.
Then $L_{n}(x) \sim (1, f(t)=\frac{t}{t-1})$, and $S_{2}(n,k;\bold{P})=\sum_{l=k}^{n}(-1)^{l}L(n,l)S_{2}(l,k)$. By Theorem 5.1,  we see that
\begin{align*}
S_{1}(n,k;\bold{P})&=\sum_{l=k}^{n}S_{1}(n,l)(-1)^{k}\big\langle{\frac{1}{k!}\Big(\frac{t}{1-t}\Big)^{k}|x^{l}\big\rangle}\\
&=\sum_{l=k}^{n}S_{1}(n,l)(-1)^{k}\sum_{j=k}^{l}L(j,k)\frac{1}{j!}\langle{t^{j}|x^{l}\rangle}\\
&=(-1)^{k}\sum_{l=k}^{n}S_{1}(n,l)L(l,k).
\end{align*}
Thus we have $S_{1}(n,k;\bold{P})=(-1)^{k}\sum_{l=k}^{n}S_{1}(n,l)L(l,k)$ (see [23]).  Now, from Theorem 5.3, we get
\begin{align*}
&\sum_{k=l}^{n}\sum_{m=k}^{n}\sum_{j=l}^{k}(-1)^{k-j}S_{1}(n,m)L(m,k)L(k,j)S_{2}(j,l)=\delta_{n,l},\\ 
&\sum_{k=l}^{n}\sum_{j=k}^{n}\sum_{m=l}^{k}(-1)^{j-l}L(n,j)S_{2}(j,k)S_{1}(k,m)L(m,l)=\delta_{n,l}, \\
&a_{n}=\sum_{k=0}^{n}\sum_{l=k}^{n}(-1)^{l}L(n,l)S_{2}(l,k)c_{k} \Longleftrightarrow c_{n}=\sum_{k=0}^{n}(-1)^{k}\sum_{l=k}^{n}S_{1}(n,l)L(l,k)a_{k},\\
&a_{n}=\sum_{k=n}^{m}\sum_{l=n}^{k}(-1)^{l}L(k,l)S_{2}(l,n)c_{k} \Longleftrightarrow c_{n}=\sum_{k=n}^{m}(-1)^{n}\sum_{l=n}^{k}S_{1}(k,l)L(l,n)a_{k}.
\end{align*}

(o) Let $\bold{P}=\left\{B_{n}(x)\right\}$ be the sequence of Bernoulli polynomials. Then $B_{n}(x) \sim \big(\frac{e^{t}-1}{t},t\big)$, and $S_{2}(n,k;\bold{P})=\sum_{l=k}^{n}S_{2}(l,k)\binom{n}{l}B_{n-l}$. By using \eqref{7B} and Theorem 5.1, we see that (see [23])
\begin{align*}
S_{1}(n,k;\bold{P})&=\frac{1}{k!}\Big\langle{t^{k}\Big|\frac{e^{t}-1}{t}(x)_{n}\Big\rangle}=\frac{1}{k!}\Big\langle{t^{k}\Big|\int_{x}^{x+1}(u)_{n} du\Big \rangle}\\
&=\frac{1}{k!}\Big\langle{t^{k}\Big|\sum_{l=0}^{n}\frac{S_{1}(n,l)}{l+1}\big((x+1)^{l+1}-x^{l+1}\big)\Big \rangle}\\
&=\frac{1}{k!}\sum_{l=k}^{n}\frac{S_{1}(n,l)}{l+1}(l+1)_{k}\Big\langle{1 \Big | (x+1)^{l+1-k}-x^{l+1-k} \Big\rangle}\\
&=\frac{1}{k!}\sum_{l=k}^{n}\frac{S_{1}(n,l)}{l+1}(l+1)_{k}.
\end{align*}
Now, by Theorem 5.3, we obtain
\begin{align*}
&\sum_{k=l}^{n}\sum_{m=k}^{n}\sum_{j=l}^{k}\frac{1}{k!}\binom{k}{j}\frac{(m+1)_{k}}{m+1}S_{1}(n,m)S_{2}(j,l)B_{k-j}=\delta_{n,l},\\
&\sum_{k=l}^{n}\sum_{j=k}^{n}\sum_{m=l}^{k}\frac{1}{l!}\binom{n}{j}\frac{(m+1)_{l}}{m+1}S_{2}(j,k)S_{1}(k,m)B_{n-j}=\delta_{n,l},\\
&a_{n}=\sum_{k=0}^{n}\sum_{l=k}^{n}S_{2}(l,k)\binom{n}{l}B_{n-l}c_{k} \Longleftrightarrow c_{n}=\sum_{k=0}^{n}\frac{1}{k!}\sum_{l=k}^{n}\frac{S_{1}(n,l)}{l+1}(l+1)_{k}a_{k},\\
&a_{n}=\sum_{k=n}^{m}\sum_{l=n}^{k}S_{2}(l,n)\binom{k}{l}B_{k-l}c_{k} \Longleftrightarrow c_{n}=\sum_{k=n}^{m}\frac{1}{n!}\sum_{l=n}^{k}\frac{S_{1}(k,l)}{l+1}(l+1)_{n}a_{k}.
\end{align*}

(p) Let $\bold{P}=\left\{E_{n}(x)\right\}$ be the sequence of Euler polynomials. Then $E_{n}(x) \sim \big(\frac{e^{t}+1}{2},t\big)$, and $S_{2}(n,k;\bold{P})=\sum_{l=k}^{n}\binom{n}{l}S_{2}(l,k)E_{n-l}$. By Theorem 5.1, we note that (see [23])
\begin{align*}
S_{1}(n,k;\bold{P})&=\frac{1}{2 k!}\big\langle{e^{t}+1\big|t^{k} (x)_{n} \rangle}
=\frac{1}{2 k!}\sum_{l=k}^{n}S_{1}(n,l)(l)_{k}\big\langle{e^{t}+1\big| x^{l-k}
\rangle}\\
&=\frac{1}{2 k!}\sum_{l=k}^{n}S_{1}(n,l)(l)_{k}(1+\delta_{l,k})
=\frac{1}{2 k!}\sum_{l=k}^{n}S_{1}(n,l)(l)_{k}+\frac{1}{2}S_{1}(n,k).
\end{align*}
Thus, by Theorem 5.3, we obtain
\begin{align*}
&\sum_{k=l}^{n}\sum_{m=k}^{n}\sum_{j=l}^{k}\frac{1}{2}\binom{m}{k}\binom{k}{j}S_{1}(n,m)S_{2}(j,l)E_{k-j}\\
&+\sum_{k=l}^{n}\sum_{j=l}^{k}\frac{1}{2}\binom{k}{j}S_{1}(n,k)S_{2}(j,l)E_{k-j}=\delta_{n,l},\\
&\sum_{k=l}^{n}\sum_{j=k}^{n}\sum_{m=l}^{k}\frac{1}{2}\binom{m}{l}\binom{n}{j}S_{2}(j,k)S_{1}(k,m)E_{n-j}\\
&+\sum_{k=l}^{n}\sum_{j=k}^{n}\frac{1}{2}\binom{n}{j}S_{2}(j,k)S_1(k,l)E_{n-j}=\delta_{n,l},\\
& a_{n}=\sum_{k=0}^{n}\sum_{l=k}^{n}\binom{n}{l}S_{2}(l,k)E_{n-l}c_{k} \Longleftrightarrow c_{n}=\sum_{k=0}^{n}\left\{\frac{1}{2 k!}\sum_{l=k}^{n}S_{1}(n,l)(l)_{k}+\frac{1}{2}S_{1}(n,k)\right\}a_{k},\\
&a_{n}=\sum_{k=n}^{m}\sum_{l=n}^{k}\binom{k}{l}S_{2}(l,n)E_{k-l}c_{k} \Longleftrightarrow c_{n}=\sum_{k=n}^{m}\left\{\frac{1}{2 n!}\sum_{l=n}^{k}S_{1}(k,l)(l)_{n}+\frac{1}{2}S_{1}(k,n)\right\}a_{k}.
\end{align*}

(q) Let $\bold{P}=\left\{(rx+s)_{n}\right\}$, with $ r \ne 0$. Then $(rx+s)_{n} \sim \big(e^{-\frac{s}{r} t}, e^{\frac{t}{r}}-1\big)$, and $S_{2}(n,k;\bold{P})=G(n,k;r,s)$ are the Gould-Hopper numbers defined by $(rx+s)_{n}=\sum_{k=0}^{n}G(n,k;r,s)(x)_{k}$. By \eqref{7B}, \eqref{17B} and Theorem 5.1, we see that (see [19])
\begin{align*}
S_{1}(n,k;\bold{P})&=\frac{1}{k!}\big\langle{e^{-\frac{s}{r}t}\big(e^{\frac{t}{r}}-1\big)^{k} \big | (x)_{n} \big \rangle}
=\frac{1}{k!}\Big\langle{(e^{t}-1)^{k} \Big | \Big(\frac{x-s}{r}\Big)_{n} \Big \rangle}\\
&=\sum_{l=k}^{n}S_{1}(n,l)\frac{1}{r^{l}}\frac{1}{k!}\big\langle{(e^{t}-1)^{k} |(x-s)^{l} \big\rangle}\\
&=\sum_{l=k}^{n}S_{1}(n,l)\frac{1}{r^{l}}\sum_{i=k}^{l}S_{2}(i,k)\frac{1}{i!}\big\langle{t^{i} |(x-s)^{l} \big\rangle}\\
&=\sum_{l=k}^{n}\sum_{i=k}^{l}\frac{(-s)^{l-i}}{r^{l}}\binom{l}{i}S_{1}(n,l)S_{2}(i,k).
\end{align*}
From Theorem 5.3, we obtain
\begin{align*}
&\sum_{k=l}^{n}\sum_{m=k}^{n}\sum_{i=k}^{m}\frac{(-s)^{m-i}}{r^{m}}\binom{m}{i}S_{1}(n,m)S_{2}(i,k)G(k,l;r,s)=\delta_{n,l},\\
&\sum_{k=l}^{n}\sum_{m=l}^{k}\sum_{i=l}^{m}G(n,k;r,s)\frac{(-s)^{m-i}}{r^{m}}\binom{m}{i}S_{1}(k,m)S_{2}(i,l)=\delta_{n,l},\\
&a_{n}=\sum_{k=0}^{n}G(n,k;r,s)c_{k} \Longleftrightarrow c_{n}=\sum_{k=0}^{n}\sum_{l=k}^{n}\sum_{i=k}^{l}\frac{(-s)^{l-i}}{r^{l}}\binom{l}{i}S_{1}(n,l)S_{2}(i,k)a_{k},\\
&a_{n}=\sum_{k=n}^{m}G(k,n;r,s)c_{k} \Longleftrightarrow c_{n}=\sum_{k=n}^{m}\sum_{l=n}^{k}\sum_{i=n}^{l}\frac{(-s)^{l-i}}{r^{l}}\binom{l}{i}S_{1}(k,l)S_{2}(i,n)a_{k}.
\end{align*}

(r) Let $\bold{P}=\left\{b_{n}(x)\right\}$ be the sequence of Bernoulli polynomials of the second kind. Then $b_{n}(x) \sim \big(\frac{t}{e^{t}-1}, e^{t}-1\big)$, and $S_{2}(n,k;\bold{P})=\binom{n}{k}b_{n-k}$,
where $b_{n}=b_{n}(0)$ are the Bernoulli numbers of the second. By using \eqref{12B} and Theorem 5.1 and noting $(x)_{n} \sim (1,e^{t}-1)$,  we have (see [23])
\begin{align*}
S_{1}(n,k;\bold{P})&=\frac{1}{k!}\Big\langle{\frac{t}{e^{t}-1} \Big| (e^{t}-1)^{k}(x)_{n}\Big\rangle}=\binom{n}{k}\Big\langle{\frac{t}{e^{t}-1} \Big| (x)_{n-k}\Big\rangle}\\
&=\binom{n}{k}\sum_{l=0}^{n-k}S_{1}(n-k,l)\Big\langle{\frac{t}{e^{t}-1} \Big| x^{l}\Big\rangle}=\binom{n}{k}\sum_{l=0}^{n-k}S_{1}(n-k,l)B_{l},
\end{align*}
where $B_{l}$ are the Bernoulli numbers. 
Now, by Theorem 5.3, we get
\begin{align*}
&\sum_{k=l}^{n}\sum_{m=0}^{n-k}\binom{n}{k}\binom{k}{l}S_{1}(n-k,m)B_{m}b_{k-l}=\delta_{n,l},\\
&\sum_{k=l}^{n}\sum_{m=0}^{k-l}\binom{n}{k}\binom{k}{l}S_{1}(k-l,m)B_{m}b_{n-k}=\delta_{n,l},\\
&a_{n}=\sum_{k=0}^{n}\binom{n}{k}b_{n-k}c_{k} \Longleftrightarrow c_{n}=\sum_{k=0}^{n}\binom{n}{k}\sum_{l=0}^{n-k}S_{1}(n-k,l)B_{l}a_{k},\\
&a_{n}=\sum_{k=n}^{m}\binom{k}{n}b_{k-n}c_{k} \Longleftrightarrow c_{n}=\sum_{k=n}^{m}\binom{k}{n}\sum_{l=0}^{k-n}S_{1}(k-n,l)B_{l}a_{k}.
\end{align*}

(s) Let $\bold{P}=\left\{C_{n}(x;a)\right\}$ be the sequence of Poisson-Charlier polynomials. Then $C_{n}(x;a) \sim \big(e^{a(e^{t}-1)}, a(e^{t}-1)\big)$, with $a \ne 0$, and $S_{2}(n,k;\bold{P})=\binom{n}{k}(-1)^{n-k}a^{-k}$. By Theorem 5.1, we get
\begin{align*}
S_{1}(n,k;\bold{P})&=\frac{a^{k}}{k!}\big\langle{e^{a(e^{t}-1)}\big | \big(e^{t}-1\big)^{k}(x)_{n}\rangle}
=a^{k}\binom{n}{k}\big\langle{e^{a(e^{t}-1)}\big | (x)_{n-k}\rangle}\\
&=a^{k}\binom{n}{k}\sum_{l=0}^{n-k}S_{1}(n-k,l)\mathrm{Bel}_{l}(a),
\end{align*}
where $\mathrm{Bel}_{n}(x)$ are the Bell polynomials given by $e^{x(e^{t}-1)}=\sum_{n=0}^{\infty}\mathrm{Bel}_{n}(x)\frac{t^{n}}{n!}$ (see [23]).\\
Now, by Theorem 5.3, we obtain
\begin{align*}
&\sum_{k=l}^{n}\sum_{m=0}^{n-k}(-a)^{k-l}\binom{n}{k}\binom{k}{l}S_{1}(n-k,m)\mathrm{Bel}_{m}(a)=\delta_{n,l},\\
&\sum_{k=l}^{n}\sum_{m=0}^{k-l}(-1)^{n-k}a^{l-k}\binom{n}{k}\binom{k}{l}S_{1}(k-l,m)\mathrm{Bel}_{m}(a)=\delta_{n,l},\\
&a_{n}=\sum_{k=0}^{n}\binom{n}{k}(-1)^{n-k}a^{-k}
c_{k} \Longleftrightarrow c_{n}=\sum_{k=0}^{n}a^{k}\binom{n}{k}\sum_{l=0}^{n-k}S_{1}(n-k,l)\mathrm{Bel}_{l}(a)a_{k},\\
&a_{n}=\sum_{k=n}^{m}\binom{k}{n}(-1)^{k-n}a^{-n}c _{k} \Longleftrightarrow c_{n}=\sum_{k=n}^{m}a^{n}\binom{k}{n}\sum_{l=0}^{k-n}S_{1}(k-n,l)\mathrm{Bel}_{l}(a)a_{k}.
\end{align*}

(t) Let $\bold{P}=\left\{p_{n}(x)\right\}$, with $p_{n}(x)=\sum_{k=0}^{n}B_{k}(x)B_{n-k}(x)$. This is not a Sheffer sequence. Here we would like to describe how to determine $S_{1}(n,k;\bold{P})$.
For this, we recall from [13] that
\begin{align}
p_{n}(x)=\frac{2}{n+2}\sum_{m=0}^{n-2}\binom{n+2}{m}B_{n-m}B_m(x)+(n+1)B_n(x).\label{2F}
\end{align}
On the one hand, by (o) we have
\begin{align}
(x)_{n}=\sum_{m=0}^{n}\gamma_{m}B_{m}(x), \label{3F}
\end{align}
where $\gamma_{m}=\frac{1}{m!}\sum_{l=m}^{n}\frac{S_{1}(n,l)}{l+1}(l+1)_{m}$.\\
On the other hand, by definition and from \eqref{2F} we also have
\begin{align}
(x)_{n}&=\sum_{k=0}^{n}S_{1}(n,k;\bold{P})p_{k}(x)\nonumber\\
&=\sum_{k=2}^{n}S_{1}(n,k;\bold{P})\frac{2}{k+2}\sum_{m=0}^{k-2}\binom{k+2}{m}B_{k-m}B_m(x) \nonumber\\
&\quad\quad+\sum_{k=0}^{n}S_{1}(n,k;\bold{P})(k+1)B_k(x)\label{4F}\\
&=\sum_{m=0}^{n-2}\sum_{k=m+2}^{n}\epsilon_{m,k}S_{1}(n,k;\bold{P})B_{m}(x)\nonumber\\
&\quad\quad+\sum_{m=0}^{n}(m+1)S_{1}(n,m;\bold{P})B_m(x),\nonumber
\end{align}
where $\epsilon_{m,k}=\frac{2}{k+2}\binom{k+2}{m}B_{k-m},\quad (0 \le m \le n-2, m+2 \le k \le n)$.\\
By comparing \eqref{3F} and \eqref{4F}, we have
\begin{align}
&\gamma_{m}=(m+1)S_{1}(n,m;\bold{P})+\sum_{k=m+2}^{n}\epsilon_{m,k}S_{1}(n,k;\bold{P}),
\quad(0 \le m \le n), \label{5F}
\end{align}
where we understand that the sum is zero for $m=n-1$ and $m=n$. Let $\Gamma$ be the column vector consisting of $\gamma_{0}, \gamma_{1}, \dots, \gamma_{n}$, and let $S$ be the column vector consisting of $S_{1}(n,0;\bold{P}), S_{1}(n,1;\bold{P}), \dots, S_{1}(n,n;\bold{P})$. Then,  in matrix form,  \eqref{5F} can be written as  $\Gamma=AS$, where $A$ is an $(n+1) \times (n+1)$ upper triangular matrix with diagonal entries $1,2, \dots , n+1$. Thus $S=A^{-1}\Gamma$.

\section{Conclusion}
In this paper, we introduced Stirling numbers of both kinds associated with sequences of polynomilas and studied them by using umbral calculus techniques, which was motivated by the observation that many important special numbers appear in the expansions of some polynomials in terms of falling factorials and vice versa. Our results were illustrated by twenty examples, both for the second kind and for the first kind. In all cases of these examples, some interesting inverse relations were obtained as immediate consequences of our results. The novelty of this paper is that it is the first paper which studies the Stirling numbers of both kinds associated with any sequence of polynomials in a unified and systematic way with the help of umbral calculus. \\
In a forthcoming paper, we study the central factorial numbers of the first kind and those of the second kind associated with sequences of polynomials with the help of umbral calculus techniques. \\
In general, it is one of our future projects to study special numbers and polynomials by making use of various tools, including combinatorial methods, generating functions, umbral calculus, $p$-adic analysis, probability theory, differential equations, analytic number theory, operator theory, special functions and so on.


\end{document}